\documentclass[letterpaper,12pt,leqno]{article}

\usepackage{fullpage}

\usepackage{hyperref}

\usepackage{amsmath}
\usepackage{amsfonts}
\usepackage{amssymb}
\usepackage{tabularx}
\usepackage{color}
\usepackage{fancyhdr}
\usepackage{tikz-cd}
\usepackage{enumitem}
\usetikzlibrary{calc}

\def\theoremparent{section}
\usepackage{notes}

\numberwithin{equation}{theorem}
\def\cA{\mathcal{A}}

\def\cZ{\mathcal{Z}}

\def\coh{\mathop{\mathrm{coh}}\nolimits}

\def\d{\partial}
\def\D{\mathrm{D}}

\def\Der{\mathrm{Der}}

\def\diag{\mathrm{diag}}

\def\gr{\mathop{\mathrm{gr}}\nolimits}

\def\k{\mathbf{k}}

\def\op{\mathrm{op}}

\def\Perf{\mathrm{Perf}}

\def\proof{\noindent{\em Proof:}\ }

\def\qed{\hfill\lower 1em\hbox{$\square$}\vskip 1em}

\def\toto{\text{\ \raise0.2em\hbox to 0pt {$\to$}\lower0.2em\hbox{$\to$}\ }}

\def\Vect{\mathrm{Vect}}

\begin{document}

\title{Tensor-Hochschild complex}
\author{Slava Pimenov, Angel Toledo}
\date{\today}
\titlepage
\maketitle

\tableofcontents
\vskip 5em

\setcounter{section}{-1}
\section{Introduction}
For an associative algebra $A$ an infinitesimal deformation of the multiplication $\mu\from A \tensor A \to A$ is given by
$$
\mu'(x, y) = \mu(x, y) + \delta(x, y) \epsilon,
$$
where $\delta$ is a also an element of $\Hom(A \tensor A, A)$. The associativity of the new multiplication $\mu'$, ignoring
terms in $\epsilon$ of degrees at least two, is linear in $\delta$ and can be seen to be equivalent to the Hochschild 2-cocycle condition
on $\delta$.

This simple observation can be upgraded to a more precise statement, the Hochschild cochain complex $C^\bullet(A, A)$ of the associative algebra $A$ possesses a $1$-shifted Lie bracket $[-,-]_G$, known
as the Gerstenhaber bracket (\cite{Gerstenhaber}). This dg-Lie algebra describes a formal deformation functor of $A$, which sends an Artinian algebra $R$
to the set of algebra maps from $R$ to the Chevalley-Eilenberg complex of $C^\bullet(A, A)[1]$ equipped with the Gerstenhaber bracket.

This can be generalized to the setting of dg-algebras, in which case the Hochschild complex controls the deformation of $A$ in the category of
curved $A_\infty$-algebras. Furthermore, it can be extended to dg-categories, as has been shown in \cite{Keller}. If a dg-category $\C$ satisfies certain
boundedness assumption on the cohomology then its Hochschild complex controls the Morita deformations of $\C$.

\begin{nparagraph}
Now, let $(\C, \tensor)$ be a monoidal category, and denote by $T^n\from \C^{\tensor n} \to \C$ for $n \ge 1$ the iterated tensor product functor. In case when
the monoidal structure is unital we also put $T^0\from * \to C$ the functor from the point category (with a single object with the trivial automorphism group) to $\C$
that sends the object to the unit $I$ of the monoidal structure. Davydov, Crane-Yetter, and Yetter (\cite{Davydov}, \cite{Yetter},\cite{Yetter2}) have constructed a complex out of endomorphisms of these functors by putting
$C^n_{DY}(\C, \tensor) = \End(T^n)$, with the differentials mimicking the differentials of the Hochschild complex. They have shown that the second cohomology group
of this complex controls the infinitesimal deformation of the identity $\Id\from \C \to \C$ as the monoidal functor, and the third cohomology controls the
deformation of the associator of the monoidal structure itself. More precisely, the complex $C^\bullet_{DY}(\C)$ is equipped with a 2-shifted dg-Lie algebra structure
describing a formal deformation functor for the associator of monoidal structure on $\C$.
\end{nparagraph}

\begin{nparagraph}
This complex also admits a more conceptual interpretation (\cite{Gainutdinov}). Let $\cZ = \cZ(\C, \tensor)$ be the Drinfeld center of the monoidal category $\C$. Objects of $\cZ$ are pairs
$(X, \rho_X)$, consisting of an object $X \in \C$ and a natural transformation expressing centrality of $X$:
$$
\rho_X\from X \tensor (-) \to (-) \tensor X.
$$
The forgetful functor $U\from\cZ \to \C$ has a left adjoint ``centralization'' functor $Z \from \C \to \cZ$. Moreover $U$ is monadic, which means that the Drinfeld center is equivalent
to the category of $UZ$-modules in $\C$. Iterating the composition of functors $U$ and $Z$ we obtain a standard resolution $P^\bullet(X)$ of an object $X \in \cZ$ by free $UZ$-modules, by setting
$P^n(X) = (ZU)^n (X)$.

The unit object $I$ of the monoidal structure is the image of object $(I, \rho_I) \in \cZ$ under forgetful functor, where $\rho_I$ is expressed as composition of the unit identities:
$$
\begin{tikzcd}
\rho_I\from I \tensor Y \ar[r, "\isom"] & Y & Y \tensor I \ar[l, "\isom"'].
\end{tikzcd}
$$
Then the Davydov-Yetter complex $C^\bullet_{DY}(\C, \tensor)$ described above is isomorphic to the complex $\Hom_\cZ(P^\bullet(I), I)$.
\end{nparagraph}

\begin{nparagraph}
The shortcoming of the Davydov-Yetter complex as a deformation complex of monoidal category $(\C, \tensor)$ is twofold. First, the underlying category $\C$ remains intact under
such deformations, and secondly even if we are only interested in the deformations of the monoidal structure leaving $\C$ intact, the complex only captures the deformation of the
associator maps, and not the deformation of tensor product on objects or maps.

The main result of this paper is the construction of the $\tensor$-Hochschild complex $TC^\bullet(\C, \tensor)$ that captures deformation of all the data of the monoidal structure
on $\C$ together with the underlying category itself. A cochain of this complex of bidegree $(p,q)$ takes a rectangle of composable arrows of height $(p+1)$ and length $q$ as in the picture
below and sends it to a map from the tensor product of starting objects of each row to the tensor product of terminal objects.

$$
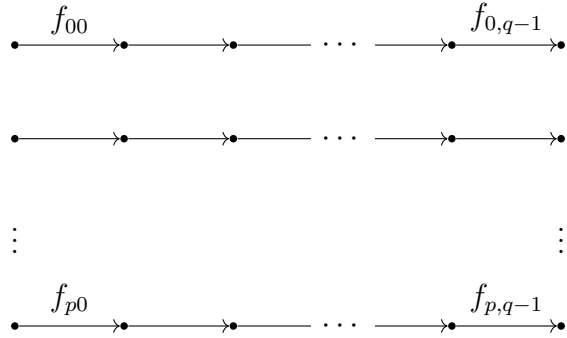
\begin{tikzpicture}
\def\u{3em}
\def\uh{3.5em}
\node (v00) at (0, 0) [circle,fill,inner sep=1pt] {};
\node (v01) at (1*\uh, 0) [circle,fill,inner sep=1pt] {};
\node (v02) at (2*\uh, 0) [circle,fill,inner sep=1pt] {};
\node (v03) at (4*\uh, 0) [circle,fill,inner sep=1pt] {};
\node (v04) at (5*\uh, 0) [circle,fill,inner sep=1pt] {};
%\node at ($(v04) + (0.5em,0)$) [right] {$P \in \cA(A_{i-1})$};

\node (v10) at (0, -1*\u) [circle,fill,inner sep=1pt] {};
\node (v11) at (1*\uh, -1*\u) [circle,fill,inner sep=1pt] {};
\node (v12) at (2*\uh, -1*\u) [circle,fill,inner sep=1pt] {};
\node (v13) at (4*\uh, -1*\u) [circle,fill,inner sep=1pt] {};
\node (v14) at (5*\uh, -1*\u) [circle,fill,inner sep=1pt] {};

\node (v20) at (0, -3*\u) [circle,fill,inner sep=1pt] {};
%\node at ($(v20) - (0,1.5em)$) {$I_{q-i+2}$};
\node (v21) at (1*\uh, -3*\u) [circle,fill,inner sep=1pt] {};
\node (v22) at (2*\uh, -3*\u) [circle,fill,inner sep=1pt] {};
\node (v23) at (4*\uh, -3*\u) [circle,fill,inner sep=1pt] {};
\node (v24) at (5*\uh, -3*\u) [circle,fill,inner sep=1pt] {};

\draw [->] (v00) -- (v01) node[midway,above] {$f_{00}$};
\draw [->] (v01) -- (v02);
\draw [->] (v02) -- (v03) node[fill=white,midway] {$\cdots$};
\draw [->] (v03) -- (v04) node[midway,above] {$f_{0,q-1}$};

\draw [->] (v10) -- (v11);
\draw [->] (v11) -- (v12);
\draw [->] (v12) -- (v13) node[fill=white,midway] {$\cdots$};
\draw [->] (v13) -- (v14);

\draw [->] (v20) -- (v21) node[midway,above] {$f_{p0}$};
\draw [->] (v21) -- (v22);
\draw [->] (v22) -- (v23) node[fill=white,midway] {$\cdots$};
\draw [->] (v23) -- (v24) node[midway,above] {$f_{p,q-1}$};

\node at (0, -2*\u) {$\vdots$};
\node at (5*\uh, -2*\u) {$\vdots$};
\end{tikzpicture}
$$

In the case when associator maps of the monoidal structure are all identities the complex $TC^\bullet$ is the totalization of a bicomplex
with vertical differential $d_0$ given by the Hochschild differential and horizontal differential $d_1$ by the differential of the Davydov-Yetter complex.
In general however, this is not a bicomplex, and the sum $d_0 + d_1$ has to be adjusted by higher terms $d_n$, for $n \ge 2$ expressing coherence of the tensor product and
the categorical composition to ensure that the total differential squares to zero. We give an explicit description of these higher terms of the differential
via certain ``admissible'' paths in Stasheff associahedra.

\end{nparagraph}

\begin{nparagraph}
The complex $TC^\bullet(\C)$ comes with natural maps relating it to the Hochschild and Davydov-Yetter complexes.
$$
\begin{tikzcd}
C^\bullet_{DY}(\C) \ar[r, into] & TC^\bullet(\C) \ar[r, epi] & C^\bullet(\C).
\end{tikzcd}
$$
We show that if the underlying category $\C$ is semisimple, then the first inclusion is a quasi-isomorphism. In other words, the usual Davydov-Yetter complex
is fully sufficient to study deformations of semisimple monoidal categories.
\end{nparagraph}

\begin{nparagraph}
A similar idea of combining Hochschild and Davydov-Yetter complexes to describe the deformations of $(\C, \tensor)$ together with the deformations
of the underlying category $\C$ has been already pursued in the literature in \cite{Shrestha}, and more recently including the dg-category case in \cite{Panero} where they construct a slightly different complex. Translating into our
language they consider a subcomplex of $TC^\bullet(\C)$ consisting of cochains, such that the two differentials $d_0$ and $d_1$ form a bicomplex.
In particular all the higher terms $d_n$ vanish on such cochains, and they do not appear in their discussion at all. While the two complexes coincide
when the associators are all identities, it seems unlikely that they are quasi-isomorphic in general.
\end{nparagraph}

\begin{nparagraph}
In the second part of the paper we consider several special cases and examples. First we study two cases where the associator maps are given by identities
and show that our complex is quasi-isomorphic to the previously constructed deformations complexes in those settings.

First, we consider the case when the monoidal category $\C$ consist of a single object $I$, the unit of the monoidal structure. The endomorphisms of this
object is a commutative algebra $A$, or more generally an $E_2$-algebra if we are working with a dg-category $\C$. The deformation complex of $A$ as an
$E_2$-algebra is given by an analog of the Hochschild complex for operad $E_2$. By comparing two spectral sequences, one converging to the operadic cohomology of $A$
and another converging to cohomology of $\tensor$-Hochschild complex of $\C$, we show that after shift of grading by one the two deformation complexes are quasi-isomorphic
(proposition \ref{prop_E2_alg}).

The other case is the monoidal category $\C$ of representations of a bialgebra $B$. The deformation of $B$ as a bialgebra (or more precisely a quasi-bialgebra) is
controlled by the Gerstenhaber-Schack complex. And we show that $TC^\bullet(\C)$ is quasi-isomorphic to this complex (section \ref{sec_bialgebras}).

The equivalence of these two special cases has been studies previously (for instance \cite{Ginot}) and our discussion here is just a reiteration of those results.
\end{nparagraph}

\begin{nparagraph}
Since the  complex $TC^\bullet(\C)$ generalizes the $E_2$-algebra cohomology complex, it is reasonable to expect a version of Deligne's conjecture to hold. Namely,
that there is an action of the chain complex of an $E_3$-operad on $TC^\bullet(\C)$. It would be interesting to make it explicit.

%It is clear that $TC^\bullet(\C)$
%is not a dg-Lie algebra, since the axioms of a monoidal category are not quadratic, but it should have a $\Lie_\infty$-algebra structure instead. The Maurer-Cartan elements
%for this structure should be thought of as describing an ``infinity-monoidal'' $A_\infty$-category. At the moment, we do not have an explicit description of the
%$\Lie_\infty$ operations of $TC^\bullet(\C)$.

In general the deformation theory of an algebraic structure is controlled by a Lie-infinity algebra (cf. \cite{KoSo}), which is a graded vector space, equipped
with a collection of operations $l_n$ of arity $n$ and degree $(2-n)$ for $n \ge 1$, satisfying certain compatibility conditions. One should think of the operation $l_1$ as the differential,
$l_2$ as the Lie bracket, $l_3$ controls the failure of the Jacobi identity to hold, etc.

Of course such Lie-infinity algebra is defined only
up to Lie-infinity quasi-isomorphism, and any such algebra can be strictified to a dg-Lie algebra, where all operations $l_n$ for $n \ge 3$ vanish. Unfortunately we do not possess the
general theory, similar to that of the algebraic operads, which would be applicable to monoidal dg-categories, and produce in a systematic way
a dg-Lie algebra controlling its deformations. Instead we are trying to make an educated guess as to what such complex should be, and it is unlikely
we will be luck enough to guess the strictified dg-Lie algebra. In fact it is not difficult to see that the complex proposed in our work doesn't admit a
dg-Lie algebra structure describing deformation of monoidal categories, therefore we expect that it is instead a Lie-infinity algebra, with non-vanishing
higher operations.

Let $L$ be a Lie-infinity algebra controlling deformations of some algebraic structure, then the collection of such structures is given by the solutions $\phi$ of the generalized Maurer-Cartan equation.
$$
\sum_{i=1}^\infty {1 \over i!} l_i(\phi, \ldots, \phi) = 0.
$$
Let $\phi$ be such a solution, then we can twist the Lie-infinity structure on $L$ to produce the tangent complex $L^\phi$ at $\phi$. Explicitly the twisted operations are given by the following
formulas (for the proof we refer to \cite{ChuLaz}).

$$
l_n^\phi(x_1, x_2, \ldots, x_n) = \sum_{i = 0}^\infty {1 \over i!} l_{n+i}(\phi, \ldots, \phi, x_1, \ldots, x_n).
$$

Our main result is an explicit description of operation $l_1^\phi$ (in other words the differential) in this twisted Lie-infinity algebra. We would like to point out here, that
we do not possess explicit descriptions of operations $l_n$ or $l_n^\phi$ in general.

In principle such deformation procedure leads outside of the class of structures considered originally, for instance by deforming an associative algebra one ends up with
a notion of an $A_\infty$ algebra. Similarly, in this case the generalized Maurer-Cartan elements for this Lie-infinity algebra should be thought of as describing an ``infinity-monoidal'' $A_\infty$-category.

Using this terminology, we would like to compare our result with the work of Shrestha (\cite{Shrestha}), where he studies the deformation theory of $R$-linear monoidal
categories over a commutative algebra $R$, in other words with Hom-spaces being $R$-modules with $R$-linear compositions.

First of all, restricting our construction of the reduced $\tensor$-Hochschild cochains $C^{pq}$ to a monoidal dg-category with Hom-spaces concentrated in degree $0$
we exactly recover the bigraded $R$-module $C^{(p,q)}$ of normalized Hochschild cochains considered by Shrestha.

Furthermore, the four differentials defined by Shrestha before theorem 2.0.38 (\cite{Shrestha}) coincide with our differential in low degrees, written out explicitly in example \ref{exa_diff}.
Thus, our result (proposition \ref{prop_complex}) provides a solution to his conjecture 2.0.39.

Next, he considers obstructions $\omega^{(n)}$ to lifting an infinitesimal deformation $\phi_1$ to a formal one
$$
\phi = \phi_0 + \phi_1 \epsilon + \phi_2 \epsilon^2 + \ldots
$$
In terms of the Lie-infinity algebra structure they are expressed as
$$
\omega^{(n)} = l_1^{\phi_0}(\phi_{n+1}) + \sum_{i_1 + i_2 = n + 1} {1\over 2!} l_2^{\phi_0}(\phi_{i_1}, \phi_{i_2}) + \cdots +
\sum_{i_1 + \cdots + i_k = n + 1} {1\over k!} l_k^{\phi_0}(\phi_{i_1}, \ldots, \phi_{i_k}) + \cdots .
$$

In our paper we do not consider the problem of obstructions, since, as these formulas show, it would require much more detailed information about the
higher Lie-infinity operations, which our combinatorial techniques do not handle. However, since we have a description of all higher terms of the differential in the twisted complex
$l_1^\phi$, one could in principle manually verify that the obstructions $\omega^{(1)}$ as was considered by Shrestha are indeed cocycles, while deforming all parts of the monoidal structure,
i.e. composition, tensor product of maps and associator, at the same time. We do not perform this calculation in our paper.
\end{nparagraph}

\begin{nparagraph}
In section \ref{sec_smooth_scheme} we consider the case of a smooth scheme $X$ and the monoidal category $\C = D^b(\coh X)$ equipped with the tensor product of coherent sheaves. We obtain a spectral sequence
relating the Hochschild cohomology of $X$ and the cohomology of $TC^\bullet(\C)$. In particular in the case when $\dim X \le 2$ this spectral sequence degenerates and
we have (theorem \ref{thm_smooth_scheme})
$$
TH^\bullet(\C) \ \isom\ \bigoplus_{p \ge 0} H^\bullet(X, S^{p+1} T_X).
$$
This also implies that all infinitesimal automorphisms of the underlying category $\C$ preserve the monoidal structure and any infinitesimal deformation of $\C$ lifts
to a deformation of $\C$ together with the monoidal structure (corollary \ref{cor_sm_scheme}).
\end{nparagraph}

\begin{nparagraph}
Finally, in the section \ref{sec_quiver} we work out a specific example. Let $Q_n$ the a quiver with two vertices and $n$ arrows going from the first vertex to the second.
We consider the bounded derived category $\C$ of representations of $Q_n$ and show that it is rigid, in the sense that cohomology of $TC^\bullet(\C)$ vanishes in all degrees.
It would be interesting to investigate the question of rigidity for other quivers.

In the particular case of $Q_2$ the underlying category $\C$ is equivalent to the category $D^b(\coh\P^1)$, but not monoidal equivalent as their tensor products differ. By comparing the deformation
complexes for these two structures we see that they lie on two different orbits of the $\Aut(\C)$ action on the space of monoidal structures.

\end{nparagraph}

\vskip 1em
The authors would like to express gratitude to the Beijing Institute of Mathematical Sciences and Applications (BIMSA) for the excellent working conditions under which this work was completed.

\vskip 5em
\section{$\tensor$-Hochschild complex of a monoidal dg-category}

Let $\C$ be a dg-category over a commutative ring $R$, so that the spaces of homomorphisms $\Hom_\C(X, Y)$ are complexes of $R$-modules.
Furthermore, assume that $\C$ is equipped with a monoidal, not necessarily unital structure $\tensor:\C \times \C \to \C$.
Specifically, we have the natural associator maps
\begin{equation}
\label{eq_assoc}
a = a_{X,Y,Z}\from (X \tensor Y) \tensor Z \to X \tensor (Y \tensor Z),
\end{equation}
and we assume for simplicity that the pentagon diagram is strictly commutative (not just up to homotopy).
\begin{equation}
\label{eq_pentagon}
\begin{tikzpicture}[commutative diagrams/every diagram]
\def\u{6em}
\node (z0) at (90+72:\u) {$((X \tensor Y) \tensor Z) \tensor W$};
\node (z1) at (90:\u) {$(X \tensor Y) \tensor (Z \tensor W)$};
\node (z2) at (90-72:\u) {$X \tensor (Y \tensor (Z \tensor W))$};
\node (z3) at (90+2*72:\u-1em) {\makebox[2em][r]{$(X \tensor (Y \tensor Z)) \tensor W$}};
\node (z4) at (90+3*72:\u-1em) {\makebox[2em][l]{$X \tensor ((Y \tensor Z) \tensor W).$}};

\path[commutative diagrams/.cd, every arrow, every label]
(z0) edge (z1)
(z1) edge (z2)
(z0) edge (z3)
(z3) edge (z4)
(z4) edge (z2);
\end{tikzpicture}
\end{equation}
Although constructions here can be generalized to the case of an ``infinity-monoidal'' $A_\infty$-category in a straightforward manner, we do not give precise definition of this notion or construction
of the complex.

\vskip 5em
\subsection{Construction of the complex.}

For $p, q \ge 0$ consider complexes of (non-reduced) $\tensor$-Hochschild cochains $\wtilde C^{pq} = \wtilde C^{pq}(\C, \tensor)$ defined as
$$
\wtilde C^{pq} = \prod_{X_{pq}} \Hom_R\left(\bigotimes_{0 \le i \le p \atop 0 \le j < q} \Hom_\C(X_{ij}, X_{i,j+1}), \Hom_\C\left(\bigotimes X_{p0}, \bigotimes X_{pq}\right) \right).
$$
Here the $\Hom_R$ and the first tensor product are the internal $\Hom$-complex and the tensor product in the category of complexes of $R$-modules respectively.

%Here the product is taken over the set of objects $X_{pq}$ of $\C$, the $\Hom_R$ and the first tensor product are the internal $\Hom$-complex and the tensor product in the category of complexes of $R$-modules respectively, and the second tensor product is taken in the category $\C$.

In what follows we will only work with the reduced cochains $C^{pq} \into \wtilde C^{pq}$ which vanish on the collections of maps $(f_{ij}\from X_{ij} \to X_{i,j+1})$ containing a column of identities:
$f_{ij} = \Id_{X_{ij}}$ for some $0 \le j < q$ and all $0 \le i \le p$.

In order to define the differential in the total space
$$
TC^n := (\Tot C^{\bullet\bullet})^n = \bigoplus_{p + q + i = n} (C^{pq})^i,
$$
we will first introduce some terminology.

\begin{nparagraph}
Let $(S, \le)$ be a partially ordered set (poset for short). A \textit{path} in $S$ is a subset $P \subset S$ totally ordered with respect to the induced order.
An \textit{aposet} $S = (S, \le, \cA)$ is a poset $(S, \le)$ and a collection $\cA$ of paths in $S$ that will be called \textit{admissible} paths.

A product of two aposets $(S, \le_S, \cA_S)$ and $(T, \le_T, \cA_T)$ is the product of posets $S \times T$, and a path $P = (p_0, \ldots, p_n) \subset S \times T$ is admissible if
and only if there exist
admissible paths $Q \subset S$ and $R \subset T$, such that every pair of consecutive elements $(p_i, p_{i+1})$ in $P$ is either of the form $((q_j, r), (q_{j+1}, r))$ for some $r \in R$ and $q_j, q_{j+1} \in Q$,
in which case we will call it \textit{horizontal}, or of the form $((q, t_j), (q, t_{j+1}))$ for some $q \in Q$ and $t_j, t_{j+1}\in R$, in which case we will call it \textit{vertical}.

We will be particularly interested in the following aposets. Let $I_n$ be the totally ordered set of $(n + 1)$ elements, with the only admissible path containing the entire set $I_n$.
Geometrically, $I_n$ corresponds to the standard $n$-dimensional simplex with the order on the vertices induced by the orientation of edges.

Next, let $A_n$ be the Stasheff associahedron of dimension $n$ considered as an oriented cellular complex, so for example $A_1$ is the segment in (\ref{eq_assoc}), $A_2$ is the pentagon in
(\ref{eq_pentagon}), and so forth. Abusing notation we will also denote by $A_n$ the set of vertices of the associahedron partially ordered by the orientation of one-dimensional edges.
We will define the set of admissible paths in $A_n$ inductively. First, we put the only admissible path in $A_0$ to be the entire set $A_0$ (so $A_0$ is isomorphic to $I_0$).

Assume we defined aposets $A_k$ for all $k < n$. The faces of associahedron $A_n$ are of the form $A_i \times A_j$ with $i + j = n - 1$ (here we put $A_0$ to be a single point) and we
equip each face with admissible paths of the product of aposets. A path $P = (p_1, \ldots, p_n) \subset A_n$ is admissible in $A_n$ if and only if $p_n$ is the terminal vertex of $A_n$ and
$(p_1, \ldots, p_{n-1})$ is an admissible path in one of the faces of $A_n$ (note that admissible paths of the faces themselves are \textit{not} admissible in $A_n$).

We would like to point out that both for $I_n$ and $A_n$ all admissible paths have the same length equal to the dimension of the corresponding cellular complex.

To illustrate the definition, $A_2$ has three admissible paths: $(v_0, v_1, v_4)$, $(v_2, v_3, v_4)$ and $(v_0, v_2, v_4)$.

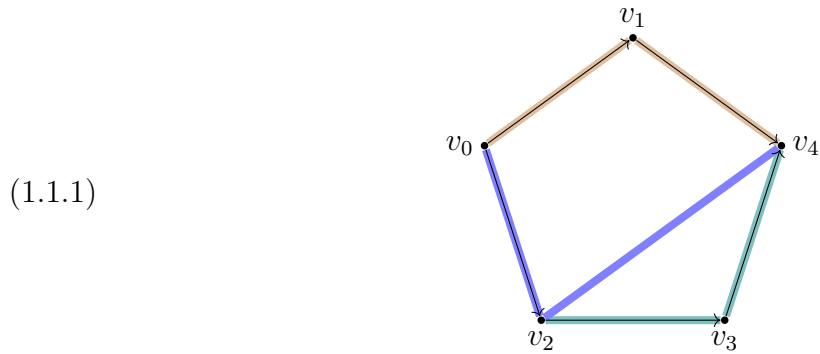
\begin{equation}
\label{equ_pentagon}
\begin{tikzpicture}[baseline=-0.25em]
\def\u{5em}
\node (v0) at (90+72:\u) [circle,fill,inner sep=1pt] {};
\node (v0n) at (90+72:\u) [left] {$v_0$};
\node (v1) at (90:\u) [circle,fill,inner sep=1pt] {};
\node (v1n) at (90:\u) [above] {$v_1$};
\node (v2) at (90+2*72:\u) [circle,fill,inner sep=1pt] {};
\node (v2n) at (90+2*72:\u) [below] {$v_2$};
\node (v3) at (90+3*72:\u) [circle,fill,inner sep=1pt] {};
\node (v3n) at (90+3*72:\u) [below] {$v_3$};
\node (v4) at (90-72:\u) [circle,fill,inner sep=1pt] {};
\node (v4n) at (90-72:\u) [right] {$v_4$};

\draw [-,color=brown!50,line width=3pt] (v0) -- (v1) -- (v4);
\draw [-,color=teal!50,line width=3pt] (v2) -- (v3) -- (v4);
\draw [-,color=blue!50,line width=3pt] (v0) -- (v2) -- (v4);

\draw [->] (v0) -- (v1);
\draw [->] (v1) -- (v4);
\draw [->] (v0) -- (v2);
\draw [->] (v2) -- (v3);
\draw [->] (v3) -- (v4);
\end{tikzpicture}
\end{equation}

\end{nparagraph}

\vskip 1em
\begin{nparagraph}
We define the differential in $TC^\bullet$ as the sum
$$
d = \delta + \sum_{i \ge 0} d_i,
$$
where $\delta$ is the internal differential of complexes $C^{pq}$ and $d_i\from C^{pq} \to C^{p + i, q - i + 1}$ are defined as follows. The component $d_0$ is the usual Hochschild differential which for a collection of maps $\{f_{ij}\}$ with $f_{ij}:X_{ij}\to X_{i,j+1}$, and $\phi\in C^{pq}$  is the morphism:
\begin{multline*}
(d_0 \phi)(\{f_{ij}\}) = (\tensor_i f_{i0}) \circ \phi(\{f_{ij}, j \neq 0\}) +\\ \sum_j (-1)^{j+1} \phi(\ldots f_{i,j+1}f_{ij}, \ldots ) + (-1)^{q+1} \phi(\{f_{ij}, j \neq q\}) \circ (\tensor_i f_{iq}).
\end{multline*}

The component $d_1$ is the Davydov-Yetter differential:
\begin{multline*}
(d_1 \phi)(\{f_{ij}\}) = (f_{0,q-1} \cdots f_{0,0}) \tensor \phi(\{f_{ij}, i \neq 0\}) +\\
\sum_i (-1)^{i+1} \phi(\ldots f_{ij} \tensor f_{i+1,j}, \ldots) + (-1)^{p+1} \phi(\{f_{ij}, i \neq p\}) \tensor (f_{p,q-1} \cdots f_{p,0}).
\end{multline*}

The remaining components express compatibility between composition and the tensor product structures. Let us look at the map $d_i\from C^{p-i, q+i-1} \to C^{pq}$ with $i \ge 2$, we need to define
$d_i \phi$ on the collection of maps $(f_{kl}\from X_{kl} \to X_{k,l+1})$, with $0 \le k \le p$ and $0 \le l < q$. Fix $j \le (p - i)$ and form a diagram $\Delta_j$ in $\C$ of shape $A_{i-1} \times I_q$,
where we consider the poset as a category in the obvious way. In the vertex $(t, n) \in A_{i-1} \times I_q$ we place tensor product
$$
t(X_{jn}, \ldots, X_{j+i,n})
$$
taken in the order prescribed by the binary tree $t$. The horizontal edges of $\Delta_j$ are given by associator maps, and the vertical edges by the suitably ordered tensor products of maps
$t(f_{jn}, \ldots, f_{j+i, n})$.

Any path in $A_{i-1} \times I_q$ determines a composable sequence of maps in $\Delta_j$. Furthermore, for an admissible path $P \subset A_{i-1} \times I_q$ we construct a collection
of maps $P(f)$ of size $(p - i) \times (q + i - 1)$. In the case of a vertical edge $(P_l, P_{l+1}) = ((t,m), (t,m+1))$ for some $t \in A_{i-1}$ we put
$$
P(f)_{kl} = \begin{cases}
f_{km},&\text{if $k < j$,}\\
t(f_{jm}, \ldots, f_{j+i, m}),&\text{if $k = j$,}\\
f_{k-i,m},&\text{if $k > j + i$.}
\end{cases}
$$
And in the case of a horizontal edge $(P_l, P_{l+1}) = ((s,m), (t,m))$ for some $m \in I_q$ we put
$$
P(f)_{kl} = \begin{cases}
\Id_{X_{km}},&\text{if $k < j$,}\\
s(X_{jm}, \ldots X_{j+i,m}) \to t(X_{jm}, \ldots X_{j+i,m}),&\text{if $k = j$,}\\
\Id_{X_{k-i,m}},&\text{if $k > j + i$.}
\end{cases}
$$

\vskip 1em

Now we define the cochain $d_i\phi$ by putting
$$
d_i\phi(\{f_{ij}\}) = \sum_j \sum_{P \in \cA(A_{i-1} \times I_q)} \pm \phi(P(f)).
$$

\end{nparagraph}

\begin{example}
\label{exa_diff}
To illustrate the definition let us provide a few explicit formulas in low degrees.
\begin{enumerate}[label=\alph*)]
\item For $d_2\from C^{01} \to C^{20}$ we have
$$
(d_2 \phi)_{X,Y,Z} = a_{X,Y,Z}^{-1} \phi(a_{X,Y,Z}) \in \End_\C((X \tensor Y) \tensor Z).
$$

\item For $d_2\from C^{02} \to C^{21}$ we have
$$
(d_2\phi)(f, g, h) = a^{-1} \phi([a, f\tensor g \tensor h]),
$$
where $f\from X \to X'$, $g\from Y \to Y'$ and $h\from Z \to Z'$. The outside associator $a = a_{X',Y',Z'}$ and the commutator is interpreted as
$$
\phi([a, f\tensor g \tensor h]) = \phi(a_{X',Y',Z'}, (f \tensor g) \tensor h) - \phi(f \tensor (g \tensor h), a_{X,Y,Z}).
$$

\item For $d_2\from C^{11} \to C^{30}$ we have
$$
(d_2 \phi)_{X,Y,Z,W} = a^{-1} \phi(a, 1) + a^{-1} \phi(1, a),
$$
where the first term is $a_{X,Y,Z}^{-1} \phi(a_{X,Y,Z}, 1_W)$ and the second term is
$$
(a_{X,Y,Z}^{-1} \tensor 1_W) a_{X,(Y\tensor Z),W}^{-1} a_{Y,Z,W}^{-1} \phi(1_X, a_{Y,Z,W}) a_{X,(Y\tensor Z),W} (a_{X,Y,Z} \tensor 1_W).
$$

\item For $d_3\from C^{02} \to C^{30}$ the differential is expressed using the admissible paths in (\ref{equ_pentagon})
$$
(d_3 \phi)_{X,Y,Z,W} = \phi(a_{14}, a_{01}) - \phi(a_{34}, a_{23}) - \phi(a_{34} a_{23}, a_{02}).
$$
Here $a_{ij}$ is the associator between vertices $v_i$ and $v_j$. More precisely, the first term is
$$
a_{(X \tensor Y), Z, W}^{-1} a_{X,Y,(Z\tensor W)}^{-1} \phi(a_{X,Y,(Z\tensor W)}, a_{(X \tensor Y), Z, W}),
$$
the second one is
$$
a_{(X \tensor Y), Z, W}^{-1} a_{X,Y,(Z\tensor W)}^{-1} \phi(1_X \tensor a_{Y,Z,W}, a_{X,(Y \tensor Z),W}) (a_{X,Y,Z} \tensor 1_W),
$$
and the last one is
$$
a_{(X \tensor Y), Z, W}^{-1} a_{X,Y,(Z\tensor W)}^{-1} \phi( (1_X \tensor a_{Y,Z,W}) a_{X,(Y \tensor Z),W}, a_{X,Y,Z} \tensor 1_W).
$$
\end{enumerate}
\end{example}

\begin{proposition}
\label{prop_complex}
$(TC^\bullet, d)$ is a complex.
\end{proposition}
\proof
Denote by $(d^2)_i \from C^{pq} \to C^{p+i,q-i+2}$ the homogeneous components of the square of the differential. We have to show that they vanish for all $i \ge 0$.

\begin{sparagraph}
First consider $(d^2)_0 \from C^{pq} \to C^{p,q+2}$. For $p = 0$ this is just the square of the standard Hochschild differential and the statement is well known. For $p \ge 1$ the proof
is similar and boils down to the associativity of composition in $\C$.
\end{sparagraph}

\begin{sparagraph}
Next consider $(d^2)_1 \from C^{pq} \to C^{p+1,q+1}$. We have $(d^2)_1 = d_0 d_1 + d_1 d_0$ and the vanishing can be checked immediately from the definitions of components $d_0$ and $d_1$
and the relation
$$
(f_1 g_1) \tensor (f_2 g_2) = (f_1 \tensor f_2) (g_1 \tensor g_2),
$$
for any two pairs of composable maps $(f_1, g_1)$ and $(f_2, g_2)$.
\end{sparagraph}

\begin{sparagraph}
Now, let us investigate $(d^2)_2 = d_2 d_0 + d_1 d_1 + d_0 d_2 \from C^{pq} \to C^{p+2,q}$, and we first focus on the case $p = 0$. For a cochain $\phi \in C^{0q}$,
one checks directly that
$$
d_1^2 \phi(f_q, \ldots, f_1 | g_q \ldots, g_1 | h_q, \ldots, h_1) = a_q^{-1} \phi( \ldots, f_i \tensor (g_i \tensor h_i), \ldots ) a_0 - \phi( \ldots (f_i \tensor g_i) \tensor h_i, \ldots ).
$$
Here $f_i \from X_{i-1} \to X_i$, $g_i \from Y_{i-1} \to Y_i$, $h_i \from Z_{i-1} \to Z_i$ and $a_0 = a_{X_0, Y_0, Z_0}$, $a_q = a_{X_q, Y_q, Z_q}$ are the corresponding associators.
We will represent graphically the two terms on the right hand side with the blue lines in the diagram below.

$$
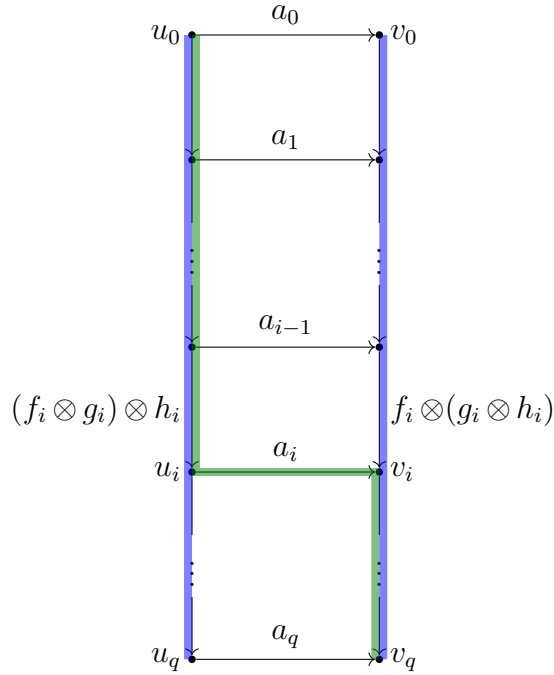
\begin{tikzpicture}
\def\u{4em}
\def\uh{6em}
\node (u0) at (0, 0) [circle,fill,inner sep=1pt] {};
\node (u0label) at (u0) [left] {$u_0$};
\node (v0) at (1*\uh, 0) [circle,fill,inner sep=1pt] {};
\node (v0label) at (v0) [right] {$v_0$};
\node (u1) at (0, -1*\u) [circle,fill,inner sep=1pt] {};
\node (v1) at (1*\uh, -1*\u) [circle,fill,inner sep=1pt] {};
\node (u2) at (0, -2.5*\u) [circle,fill,inner sep=1pt] {};
\node (v2) at (1*\uh, -2.5*\u) [circle,fill,inner sep=1pt] {};
\node (u3) at (0, -3.5*\u) [circle,fill,inner sep=1pt] {};
\node (u3label) at (u3) [left] {$u_i$};
\node (v3) at (1*\uh, -3.5*\u) [circle,fill,inner sep=1pt] {};
\node (v3label) at (v3) [right] {$v_i$};
\node (u4) at (0, -5*\u) [circle,fill,inner sep=1pt] {};
\node (u4label) at (u4) [left] {$u_q$};
\node (v4) at (1*\uh, -5*\u) [circle,fill,inner sep=1pt] {};
\node (v4label) at (v4) [right] {$v_q$};

\draw [->] (u0) -- (v0) node[above,midway] {$a_0$};
\draw [->] (u1) -- (v1) node[above,midway] {$a_1$};
\draw [->] (u2) -- (v2) node[above,midway] {$a_{i-1}$};
\draw [->] (u3) -- (v3) node[above,midway] {$a_i$};
\draw [->] (u4) -- (v4) node[above,midway] {$a_q$};

\draw [->] (u0) -- (u1);
\draw [->] (v0) -- (v1);
\draw [->] (u1) -- (u2) node[fill=white,midway] {$\vdots$};
\draw [->] (v1) -- (v2) node[fill=white,midway] {$\vdots$};
\draw [->] (u2) -- (u3) node[left,midway] {$(f_i \tensor g_i) \tensor h_i$};
\draw [->] (v2) -- (v3) node[right,midway] {$f_i \tensor (g_i \tensor h_i)$};
\draw [->] (u3) -- (u4) node[fill=white,midway] {$\vdots$};
\draw [->] (v3) -- (v4) node[fill=white,midway] {$\vdots$};

\draw [-,color=blue,opacity=0.5,line width=3pt] ($(u0) + (-1.5pt,0)$) -- ($(u4) + (-1.5pt,0)$);
\draw [-,color=green!50!black,opacity=0.5,line width=3pt] ($(u0) + (+1.5pt,0)$) -- ($(u3) + (1.5pt,0)$) -- ($(v3) + (-1.5pt,0)$) -- ($(v4) + (-1.5pt,0)$);
\draw [-,color=blue,opacity=0.5,line width=3pt] ($(v0) + (1.5pt,0)$) -- ($(v4) + (1.5pt,0)$);

\end{tikzpicture}
$$

Furthermore, every term coming from the composition $d_2 d_0$ corresponds to a path in the aposet $A_1 \times I_q$ obtained from an admissible path by removing a single vertex.
Similarly, terms in the composition $d_0 d_2$ correspond to paths obtained from an admissible path by removing a single vertex other than one of the two ``corner'' vertices ($u_i$ and $v_i$
for the green path in the picture above). Cancelling out the similar paths we find that the contribution of $d_0 d_2 + d_2 d_0$ consists of the two vertical blue paths already mentioned,
and paths of the form $(u_0, \ldots, u_{i-1}, v_i, \ldots, v_q)$. But each of the latter paths can be obtained from two different admissible paths: $(u_0, \ldots, u_i, v_i, \ldots, v_q)$
by removing vertex $u_i$ and from $(u_0, \ldots, u_{i-1}, v_{i-1}, v_i, \ldots, v_q)$ by removing $v_{i-1}$. Hence, all of them cancel out and we conclude that $(d^2)_2$ vanishes on $C^{0q}$.
\end{sparagraph}

\begin{sparagraph}
For $(d^2)_2\from C^{pq} \to C^{p+2,q}$ with $p > 0$ we once again first look at the composition $d_1^2$. A standard argument shows that terms
$$
\phi( \ldots, f_{ij} \tensor f_{i+1,j}, \ldots, f_{kj} \tensor f_{k+1,j}, \ldots )
$$
cancel each other out for disjoint pairs of indices $(i, i+1)$ and $(k, k+1)$. Similarly, terms of the following types are also all canceled out.
\begin{align*}
(f_{0,q-1} \cdots f_{0,0}) &\tensor \phi(\ldots, f_{ij} \tensor f_{i+1,j}, \ldots ),\\
(f_{0,q-1} \cdots f_{0,0}) \tensor (f_{1,q-1} \cdots f_{1,0}) &\tensor \phi(f_{2j}, \ldots ),\\
\phi(\ldots, f_{ij} \tensor f_{i+1,j}, \ldots ) &\tensor (f_{p,q-1} \cdots f_{p0}),\\
\phi(\ldots, f_{p-2,j}) &\tensor (f_{p-1,q-1} \cdots f_{p-1,0}) \tensor (f_{p,q-1} \cdots f_{p,0}).
\end{align*}

We are left with the terms
$$
\phi( \ldots, (f_{jk} \tensor f_{j+1,k}) \tensor f_{j+2,k}, \ldots ) \quad\textrm{and}\quad \phi( \ldots, f_{jk} \tensor (f_{j+1,k} \tensor f_{j+2,k}), \ldots ).
$$
Combining this with the definition of the component of the differential $d_2$ we see that $(d^2)_2$ splits into sum over $0 \le j \le (p - 2)$ and each summand can be shown
to vanish using the same argument as in step (c).
\end{sparagraph}

\begin{sparagraph}
Let us consider component $(d^2)_i\from C^{0,i-2} \to C^{i,0}$ for $i > 2$. The terms coming from the composition $d_i d_0$ have the form
$(d_0\phi)(f_{i-1}, \ldots, f_1)$, where maps $f_k$ form an admissible path in the associahedron $A_{i-1}$. In other words, these terms correspond to paths in $A_{i-1}$
obtained from an admissible path by removing one of the vertices.

Terms coming from the composition $d_1 d_{i-1}$ are of two types. First, we have
\begin{equation}
\label{prop_complex_terms}
1_{X_0} \tensor \phi(g_{i-2}, \ldots, g_1) \quad \textrm{and} \quad \phi(h_{i-2}, \ldots, h_1) \tensor 1_{X_i},
\end{equation}
where collections of maps $(g_k)$ and $(h_k)$ form an admissible path in the associahedra $A_{i-2}$ constructed respectively on $(X_1, \ldots, X_i)$ and $(X_0, \ldots, X_{i-1})$.
The other terms in $d_1 d_{i-1}$ correspond to admissible paths on the faces of $A_{i-1}$ isomorphic to $A_{i-2}$ indexed by trees $t_k$ for $1 \le k \le (i-1)$:

$$
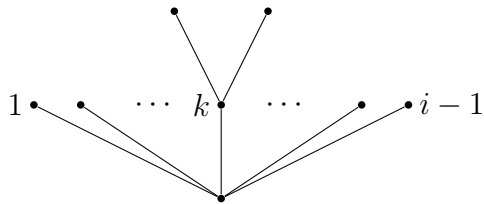
\begin{tikzpicture}
\def\u{3em}
\node (root) at (0, 0)  [circle,fill,inner sep=1pt] {};
\node (v1) at (-2 * \u, \u)  [circle,fill,inner sep=1pt] {};
\node (v1label) at (v1) [left] {$1$};
\node (v2) at (-1.5 * \u, \u)  [circle,fill,inner sep=1pt] {};
\node (v3) at (0 * \u, \u)  [circle,fill,inner sep=1pt] {};
\node (v3label) at (v3) [left] {$k$};
\node (v4) at (1.5 * \u, \u)  [circle,fill,inner sep=1pt] {};
\node (v5) at (2 * \u, \u)  [circle,fill,inner sep=1pt] {};
\node (v5label) at (v5) [right] {$i-1$};
\node (w1) at (-0.5 * \u, 2*\u)  [circle,fill,inner sep=1pt] {};
\node (w2) at (0.5 * \u, 2*\u)  [circle,fill,inner sep=1pt] {};
\node at (-0.7 * \u, \u) {$\cdots$};
\node at (0.7 * \u, \u) {$\cdots$};

\draw [-] (v1) -- (root);
\draw [-] (v2) -- (root);
\draw [-] (v3) -- (root);
\draw [-] (v4) -- (root);
\draw [-] (v5) -- (root);
\draw [-] (w1) -- (v3);
\draw [-] (w2) -- (v3);
\end{tikzpicture}
$$

For the composition $d_{i-1} d_1$ there are three type of terms. First, we have
$$
\phi(1_{X_0}, \ldots, 1_{X_0}) \tensor (g_{i-2} \cdots g_1) \quad \textrm{and} \quad (h_{i-2} \cdots h_1) \tensor \phi(1_{X_i}, \ldots, 1_{X_i}).
$$
These vanish because we are working with the reduced cochains $\phi$. Next, there are terms of the same form as in (\ref{prop_complex_terms}) which will cancel out. Finally,
there are terms corresponding to the admissible paths on the remaining two faces of $A_{i-1}$ isomorphic to $A_{i-2}$, which are indexes by the trees
$$
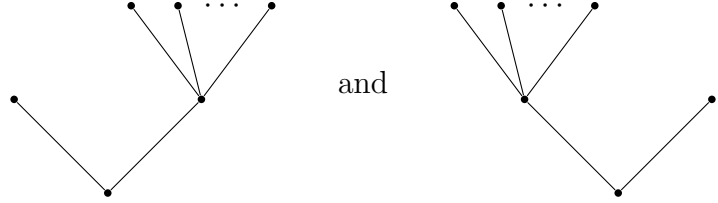
\begin{tikzpicture}[baseline=(current bounding box.east)]
\def\u{3em}
\node (root) at (0, 0)  [circle,fill,inner sep=1pt] {};
\node (w1) at (-1 * \u, \u)  [circle,fill,inner sep=1pt] {};
\node (w2) at (1 * \u, \u)  [circle,fill,inner sep=1pt] {};
\node (v1) at (0.25 * \u, 2*\u)  [circle,fill,inner sep=1pt] {};
%\node (v1label) at (v1) [left] {$1$};
\node (v2) at (0.75 * \u, 2*\u)  [circle,fill,inner sep=1pt] {};
\node (v3) at (1.75 * \u, 2*\u)  [circle,fill,inner sep=1pt] {};
%\node (v3label) at (v3) [left] {$k$};

\node at (1.25 * \u, 2*\u) {$\cdots$};

\draw [-] (v1) -- (w2);
\draw [-] (v2) -- (w2);
\draw [-] (v3) -- (w2);
\draw [-] (w1) -- (root);
\draw [-] (w2) -- (root);
\end{tikzpicture}
\quad\quad \textrm{and} \quad\quad
\begin{tikzpicture}[baseline=(current bounding box.east)]
\def\u{3em}
\node (root) at (0, 0)  [circle,fill,inner sep=1pt] {};
\node (w1) at (-1 * \u, \u)  [circle,fill,inner sep=1pt] {};
\node (w2) at (1 * \u, \u)  [circle,fill,inner sep=1pt] {};
\node (v1) at (-1.75 * \u, 2*\u)  [circle,fill,inner sep=1pt] {};
%\node (v1label) at (v1) [left] {$1$};
\node (v2) at (-1.25 * \u, 2*\u)  [circle,fill,inner sep=1pt] {};
\node (v3) at (-0.25 * \u, 2*\u)  [circle,fill,inner sep=1pt] {};
%\node (v3label) at (v3) [left] {$k$};

\node at (-0.75 * \u, 2*\u) {$\cdots$};

\draw [-] (v1) -- (w1);
\draw [-] (v2) -- (w1);
\draw [-] (v3) -- (w1);
\draw [-] (w1) -- (root);
\draw [-] (w2) -- (root);
\end{tikzpicture}
$$
\vskip 1em

Terms coming from all the other compositions $d_m d_n$, with $m, n \ge 2$ correspond to the admissible paths on the remaining faces of $A_{i-1}$. Specifically contributions from $d_m d_n$ and $d_n d_m$
cover all the faces isomorphic to $A_{m-1} \times A_{n-1}$.

For any path $P$ let us denote by $\delta(P)$ the set of paths obtained from $P$ by removing a single vertex. Also denote by $A_{i-1}(k)$ the set of codimension $k$ faces of $A_{i-1}$.
We need to show the following lemma.

\begin{lemma}
We have an equality in the free abelian group generated by paths in $A_{i-1}$
$$
\sum_{Q \in \delta(P) \atop P \in \cA(A_{i-1})} \pm Q\ \  =\ \  \sum_{Q \in \cA(F) \atop F \in A_{i-1}(1)} \pm Q.
$$
\end{lemma}
\proof
We proceed by induction on $i$. For $A_1$ the statement is obvious.

Suppose we established the statement for $A_k$ and $A_l$, then it is also true for $A_k \times A_l$. Indeed, consider an admissible path $P$ in $A_k \times A_l$ and let $v \in P$ be the removed
vertex. If $v$ sits between two horizontal or two vertical edges of $P$ then the question reduces to that of $A_k$ or $A_l$ respectively. If $v$ sits between a horizontal and a vertal edge of $P$
then the path $P - \{v\}$ can also be obtained from another admissible path $P'$ and the two terms cancel out.

We will call a codimension 1 face of associahedron $A_{i-1}$ \textit{terminal} if it contains the terminal vertex. By definition, if $P$ is an admissible path in $A_{i-1}$,
then by removing the last vertex from $P$ we obtain an admissible path on one of the non-terminal faces of $A_{i-1}$, and this will cover all admissible paths on non-terminal faces.

Let $Q$ be obtained from an admissible path $P$ by removing a non-terminal vertex, and consider $\wbar Q$ obtained from $Q$ by also removing the terminal vertex. So $\wbar Q$ can be obtained
from an admissible path of a non-terminal face by removing one vertex. Fix a non-terminal face $F$, then by inductive assumption sum over all $\wbar Q$ belonging to $F$ equals to the sum
over all admissible paths of codimension 1 faces of $F$. However, each such face belongs to two codimension 1 faces of $A_{i-1}$, say $F$ and $F'$. If both of them are non-terminal, then the path $Q$ cancels
out with another path obtained from some admissible path $P'$ in $A_{i-1}$, such that $\wbar P'$ lies on $F'$. Finally, if $F'$ is a terminal face, then $Q$ is an admissible path on $F'$.
This completes the proof of the lemma. 
\qed

\end{sparagraph}

\begin{sparagraph}
Now, let us consider the general case of $(d^2)_i\from C^{pq} \to C^{p+i,q-i+2}$. First of all, arguing as in step (d), all terms in the part of the composition $d_m d_n$ corresponding to a pair of
disjoint blocks of indices $(j, j+1, \ldots, j+m)$ and $(k, k+1, \ldots, k+n)$ cancel each other out. So the only contribution can come when the two blocks merge together into a single block
$(j, j+1, \ldots, j+m+n)$. It is enough to consider each such block separately, which reduces the problem to the case $p = 0$.

Assume from now on that $p = 0$, arguing as in step (c), terms in the composition $d_i d_0$ correspond to the paths in $\delta(P)$ for some admissible path $P$ in $A_{i-1} \times I_{q-i+2}$.
The composition $d_0 d_i$ corresponds only to those paths in $\delta(P)$ that are obtained from $P$ by removing a vertex between two vertical edges. And compositions $d_m d_n$ give all admissible
paths in $F \times I_{q-i+2}$ for all codimension 1 faces $F$ in $A_{i-1}$.
$$
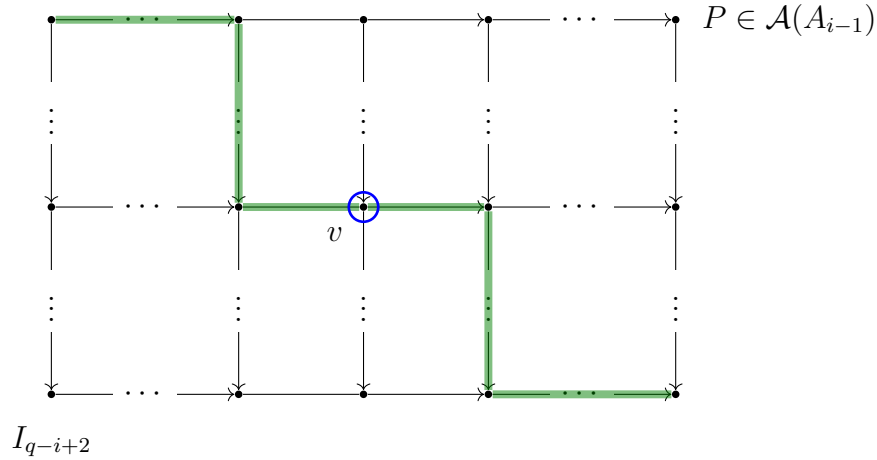
\begin{tikzpicture}
\def\u{4em}
\def\uh{4em}
\node (v00) at (0, 0) [circle,fill,inner sep=1pt] {};
\node (v01) at (1.5*\uh, 0) [circle,fill,inner sep=1pt] {};
\node (v02) at (2.5*\uh, 0) [circle,fill,inner sep=1pt] {};
\node (v03) at (3.5*\uh, 0) [circle,fill,inner sep=1pt] {};
\node (v04) at (5*\uh, 0) [circle,fill,inner sep=1pt] {};
\node at ($(v04) + (0.5em,0)$) [right] {$P \in \cA(A_{i-1})$};

\node (v10) at (0, -1.5*\u) [circle,fill,inner sep=1pt] {};
\node (v11) at (1.5*\uh, -1.5*\u) [circle,fill,inner sep=1pt] {};
\node (v12) at (2.5*\uh, -1.5*\u) [circle,fill,inner sep=1pt] {};
\node (v13) at (3.5*\uh, -1.5*\u) [circle,fill,inner sep=1pt] {};
\node (v14) at (5*\uh, -1.5*\u) [circle,fill,inner sep=1pt] {};

\node (v20) at (0, -3*\u) [circle,fill,inner sep=1pt] {};
\node at ($(v20) - (0,1.5em)$) {$I_{q-i+2}$};
\node (v21) at (1.5*\uh, -3*\u) [circle,fill,inner sep=1pt] {};
\node (v22) at (2.5*\uh, -3*\u) [circle,fill,inner sep=1pt] {};
\node (v23) at (3.5*\uh, -3*\u) [circle,fill,inner sep=1pt] {};
\node (v24) at (5*\uh, -3*\u) [circle,fill,inner sep=1pt] {};

\draw [->] (v00) -- (v01) node[fill=white,midway] {$\cdots$};
\draw [->] (v01) -- (v02) -- (v03);
\draw [->] (v03) -- (v04) node[fill=white,midway] {$\cdots$};

\draw [->] (v10) -- (v11) node[fill=white,midway] {$\cdots$};
\draw [->] (v11) -- (v12) -- (v13);
\draw [->] (v13) -- (v14) node[fill=white,midway] {$\cdots$};

\draw [->] (v20) -- (v21) node[fill=white,midway] {$\cdots$};
\draw [->] (v21) -- (v22) -- (v23);
\draw [->] (v23) -- (v24) node[fill=white,midway] {$\cdots$};

\draw [->] (v00) -- (v10) node[fill=white,midway] {$\vdots$};
\draw [->] (v10) -- (v20) node[fill=white,midway] {$\vdots$};
\draw [->] (v01) -- (v11) node[fill=white,midway] {$\vdots$};
\draw [->] (v11) -- (v21) node[fill=white,midway] {$\vdots$};
\draw [->] (v02) -- (v12) node[fill=white,midway] {$\vdots$};
\draw [->] (v12) -- (v22) node[fill=white,midway] {$\vdots$};
\draw [->] (v03) -- (v13) node[fill=white,midway] {$\vdots$};
\draw [->] (v13) -- (v23) node[fill=white,midway] {$\vdots$};
\draw [->] (v04) -- (v14) node[fill=white,midway] {$\vdots$};
\draw [->] (v14) -- (v24) node[fill=white,midway] {$\vdots$};

\draw [-,color=green!50!black,opacity=0.5,line width=3pt] (v00) -- (v01) -- (v11) -- (v12) -- (v13) -- (v23) -- (v24);
\node (v) at (v12) [circle,draw,line width=1pt,color=blue] {};
\node at (v) [below left=0.3em] {$v$};

\end{tikzpicture}
$$

Similar to step (c), the sum $d_i d_0 + d_0 d_i$ has contributions only from those paths in $\delta(P)$ which are obtained from an admissible path $P$ by removing a vertex between two horizontal edges
(the blue vertex $v$ in the picture above), the first vertex of $P$ if it starts with a horizontal edge or the last vertex of $P$ if it ends with a horizontal edge.
Finally, arguing as in step (e), all such paths cancel out with the paths from the remaining compositions $d_m d_n$ with $m, n \ge 1$. This completes the proof of the proposition.
\qed
\end{sparagraph}

\vskip 1em
\begin{nparagraph}[Relation to Hochschild and Davydov-Yetter cohomology.]
\label{par_filtrations}
First, consider the decreasing filtration of $TC^\bullet$ by slicing it in the vertical direction:
$$
F_j(TC^\bullet) = \bigoplus_{p \ge j} C^{pq}.
$$
The associated graded piece $\gr^0_F(TC^\bullet)$ is isomorphic to the relative Hochschild complex of dg-category $\C$ over $R$. In particular in the case when $R = k$ is a field it is the
usual Hochschild complex of $\C$, and when $\C$ is a category with a single object with the endomorphism ring $E$ we recover the relative Hochschild cohomology complex
of algebras $HH^\bullet(E \ |\  R)$.

The natural map of cohomologies induced by the surjective map $TC^\bullet \epi \gr^0_F(TC^\bullet)$ corresponds to forgetting the monoidal structure of $(\C, \tensor)$.

Next let us consider the increasing filtration of $TC^\bullet$ by canonically slicing it in the horizontal direction:
$$
G^j(TC^\bullet) = \bigoplus_{q < j} C^{pq} \oplus \bigoplus \Ker(d_0\from C^{pj} \to C^{p,j+1}).
$$
It is immediate to see that the vanishing of the differential $d_0\from C^{p0} \to C^{p1}$ expresses the naturality of the collection of endomorphisms in $C^{p0}$, thus
$$
G^0(TC^\bullet)^{p0} = \End_\C(\Id_\C^{\tensor p}).
$$
And the complex $G^0(TC^\bullet)$ is isomorphic to the Davydov-Yetter complex of the monoidal category $\C$. We denote $H^\bullet_{DY}(\C)$ the cohomology of this complex.

\end{nparagraph}

\begin{proposition}
Let $\C$ be a semisimple finitely generated dg-category, then the natural map
$$
H^\bullet_{DY}(\C) \to H^\bullet(TC^\bullet)
$$
induced by the inclusion $G^0(TC^\bullet) \into TC^\bullet$ is an isomorphism.
\end{proposition}
\proof
Let $X_1, \ldots, X_n$ be a collection of simple objects generating $\C$. The endomorphism algebra $A = \End_\C(\bigoplus X_i)$ is isomorphic to the direct sum $\k^n$.
Similarly, the tensor product $A^{\tensor p}$ is isomorphic to $\k^{n^p}$.

Consider the spectral sequence associated to the filtration $F$. Since the collection of $\{X_i\}$ generates $\C$, on the the layer $E_1$ of the spectral sequence the column $p$
is isomorphic to the Hochschild cohomology
$$
E_1^{pq} \ \isom\ HH^\bullet \left(A^{\tensor (p+1)}, \End_\C((\bigoplus X_i)^{\tensor (p+1)}) \right).
$$
These cohomology groups vanish for $q > 1$, and for $q = 0$ they are isomorphic to $G^0(TC^\bullet)$. The differential on $E_1$ coincides with the restriction of the
differential of $TC^\bullet$ to $G^0(TC^\bullet)$ and the statement immediately follows.
\qed

\begin{nparagraph}
It will be convenient to consider the subcomplex $F_1(TC^\bullet)$ of $TC^\bullet$. Intuitively, it controls deformations of the monoidal structure of $\C$, which leave
the underlying category $\C$ intact. To simplify notations we will write
$$
TH^n(\C) = H^n(TC^\bullet(\C)),\quad \text{and} \quad \wbar{TH}^n(\C) = H^n(F_1(TC^\bullet)).
$$
\end{nparagraph}

\begin{nparagraph}
More generally, instead of working over a commutative base ring $R$ one could construct the $\tensor$-Hochschild complex over an $E_2$-algebra base.
\end{nparagraph}

\begin{nparagraph}
Let us give interpretation of components of a 2-cocycle in $\phi \in TC^2(\C)$ that describes a deformation of a monoidal structure on $\C$. Recall that each component
$C^{pq}$ is a complex and we denote by $(C^{pq})^r$ the term of cohomological degree $r$.
\begin{enumerate}[label=\alph*)]
\item $\phi^{02}_0 \in (C^{02})^0$~--- describes the deformation of the composition of arrows.
\item $\phi^{11}_0 \in (C^{11})^0$~--- describes the deformation of the tensor product of two maps.
\item $\phi^{20}_0 \in (C^{20})^0$~--- describes the deformation of the associator.
\item $\phi^{01}_1 \in (C^{01})^1$~--- the space $(C^{01})^1$ consists of maps
$$
\phi^{01}_1(X,Y) \from \Hom_\C(X, Y) \to \Hom_\C(X, Y)[1],
$$
for any pair of objects $X, Y \in \C$.
This describes the deformation of the differential on $\Hom$-complexes of $\C$ by setting for a map $f\from X \to Y$ the new differential $d'f = df + \phi^{01}_1(X, Y)(f)$.
\item $\phi^{10}_1 \in (C^{10})^1$~--- the space $(C^{10})^1$ consists of elements $\phi^{10}_1(X,Y) \in \End^1(X \tensor Y)$, which control the infinitesimal deformation of
the object $X \tensor Y$.
\item $\phi^{00}_2 \in (C^{00})^2$~--- describes the deformation of the curvature of $\C$, which is compatible with the monoidal structure, in the sense that $c_{X \tensor Y} = 1_X \tensor c_Y + c_X \tensor 1_Y$.
\end{enumerate}

Furthermore, if $\Hom$-complexes in $\C$ have terms in negative degrees then the axioms of the deformed monoidal structure will be satisfied only up to homotopy. For instance,
\begin{enumerate}[label=\alph*)]
\item $\phi^{03}_{-1} \in (C^{03})^{-1}$~--- describes the associator for the composition of arrows in $\C$.
\item $\phi^{12}_{-1} \in (C^{12})^{-1}$~--- describes the homotopy expressing distributivity between the tensor product and the composition or arrows $(f_1 g_1) \tensor (f_2 g_2) \sim (f_1 \tensor f_2) (g_1 \tensor g_2)$.
\item $\phi^{21}_{-1} \in (C^{21})^{-1}$~--- describes the associator for the tensor product of arrows in $\C$.
\item $\phi^{30}_{-1} \in (C^{30})^{-1}$~--- describes the pentagon homotopy.
\end{enumerate}
And so forth in lower degrees.
\end{nparagraph}

\vskip 5em
\subsection{Unital version of $\tensor$-Hochschild complex.}
\label{sec_unital}
In case when the monoidal category $(\C, \tensor)$ has a unit object $I$ we can extend the $\tensor$-Hochschild complex defined above by adding one more column:
$$
C^{-1,q} = \End_\C(I)
$$
for all $q \ge 0$. We define the differential $d_0\from C^{-1,q} \to C^{-1,q+1}$ to be $0$ when $q$ is even and identity when $q$ is odd.
Following the same pattern as with the Davydov-Yetter differential $d_1$ we define $d_1\from C^{-1,q} \to C^{0,q}$ by
$$
(d_1 g)(f_0, \ldots, f_{q-1}) = g \tensor (f_{q-1} \cdots f_0) - (f_{q-1} \cdots f_0) \tensor g,
$$
for any $g \in \End_\C(I) = C^{-1,q}$.

Commutativity of the squares
$$
\begin{tikzcd}
I \tensor X \ar[r, "1_I \tensor f"] \ar[d, "g \tensor 1_X"'] & I \tensor Y \ar[d, "g \tensor 1_Y"] \\
I \tensor X \ar[r, "1_I \tensor f"] & I \tensor Y
\end{tikzcd}
$$
and the similar squares for the tensoring with $I$ from the right, immediately shows that the composition $(d_0 d_1 + d_1 d_0) \from C^{-1,q} \to C^{0,q+1}$ is zero.

Furthermore from the definition of the unital monoidal structure the commutativity of triangles

$$
\begin{tikzcd}[column sep=1em,row sep=4em]
(X \tensor I) \tensor Y \ar[rr, "a"] \ar[rd, "r_X \tensor 1_Y"'] && X \tensor (I \tensor Y) \ar[ld, "1_X \tensor l_Y"] \\
& X \tensor Y &
\end{tikzcd}
$$
as well as the following two types of triangles which can be easily deduced from the one above:
$$
\begin{tikzcd}[column sep=1em,row sep=4em]
(I \tensor X) \tensor Y \ar[rr, "a"] \ar[rd, "l_X \tensor 1_Y"'] && I \tensor (X \tensor Y) \ar[ld, "l_{X \tensor Y}"] \\
& X \tensor Y &
\end{tikzcd}
\quad\quad
\begin{tikzcd}[column sep=1em,row sep=4em]
(X \tensor Y) \tensor I \ar[rr, "a"] \ar[rd, "r_{X \tensor Y}"'] && X \tensor (Y \tensor I) \ar[ld, "1_X \tensor r_Y"] \\
& X \tensor Y &
\end{tikzcd}
$$
we also find that $d_1^2\from C^{-1,q} \to C^{1,q}$ is
\begin{multline*}
(d_1^2 h)(f,g) = ((h \tensor f) \tensor g - (f \tensor h) \tensor g) -\\ (h \tensor (f \tensor g) - (f \tensor g) \tensor h) + (f \tensor (h \tensor g) - f \tensor (g \tensor h)) = 0.
\end{multline*}
Here $f = f_{q-1} \cdots f_0$ and $g = g_{q-1} \cdots g_0$ and $h \in \End_\C(I)$.

We will denote this extended $\tensor$-Hochschild complex by $UTC^\bullet = UTC^\bullet(\C, \tensor, I)$.

\begin{nparagraph}
By taking the vertical filtration $F_p(UTC)$ as was defined in the previous section, we see that the associated graded piece $\gr^{-1}_F(UTC)$ is an acyclic complex, except in degree $0$.
Therefore, we have the long exact sequence of cohomologies
$$
\begin{tikzcd}
\cdots \ar[r] & H^n(TC) \ar[r] & H^n(UTC) \ar[r] & H^{n+1}(\End_\C(I)) \ar[r] & \cdots
\end{tikzcd}
$$

For convenience we will write
$$
UTH^n(\C) = H^n(UTC^\bullet(\C)).
$$
\end{nparagraph}

\vskip 5em
\section{Special cases and examples}

In this section we will work over a ground field $\k$ of characteristic $0$.

\subsection{Commutative algebra.}
\label{subsec_sm_comm}
Consider a monoidal dg-category $\C$ consisting of a single object $I$, which is the unit for the monoidal structure. Then the space of endomorphisms of $I$ is an $E_2$ algebra, that we will
denote by $A$. In this case, the spaces of the $\tensor$-Hochschild cochains are
$$
C^{pq} = \Hom_k\left(\left(A^{\tensor (p+1)} \right)^{\tensor q}, \Hom_\C(I^{\tensor (p+1)}, I^{\tensor (p+1)}) \right).
$$
The target of the hom-space is of course just $A$, and it is equipped with the obvious structure of $A^{\tensor (p+1)}$-bimodule.

Consider the spectral sequence associated to the stupid vertical filtration $F_\bullet(TC^\bullet)$ of the $\tensor$-Hochschild complex defined in \ref{par_filtrations}.
The differentials in $E_0^{p\bullet}$ are just the Hochschild differentials for the algebra $A^{\tensor (p+1)}$ with coefficients in the bimodule $A$. Therefore, we have
\begin{equation}
\label{equ_sp_seq_alg}
E_1^{pq} = HH^q(A^{\tensor (p+1)}, A) \Rightarrow TH^{p+q}(\C).
\end{equation}

Let us study the case when $A$ is a smooth commutative algebra of finite type.

\begin{proposition}
\label{prop_sm_comm}
If $A$ is a smooth commutative algebra, then
$$
TH^\bullet(A) \ \isom\ A \tensor \bigoplus_{p \ge 1} S_A^{p+1}(\Der(A))[-2p-1].
$$
\end{proposition}
\proof
Since $A$ is smooth, the module of differentials $\Omega^1_A$
as well as the module of derivations $\Der(A)$ are finitely generated projective $A$-modules.
Localizing if necessary we may assume them to be free, and we will write $\Der(A) \isom A \tensor V$, where $V$ is a vector space of dimension equal to the dimension of $A$.

Since $A$ is smooth commutative algebra the tensor products $A^{\tensor (p+1)}$ are also smooth commutative algebras and we have
$$
\Der(A^{\tensor (p+1)}, A) \ \isom\  A \tensor V^{\oplus (p+1)}.
$$
Hochschild-Kostant-Rosenberg theorem gives us identification
$$
HH^q(A^{\tensor (p+1)}, A) \ \isom\  \Lambda^q_A \Der(A^{\tensor (p+1)}, A) \ \isom\  A \tensor \Lambda^q(V^{\oplus (p+1)}).
$$

Let $\d x$ by an element of $V$, then we will write $\d x_i \in A \tensor \Lambda^q(V^{\oplus (p+1)})$, $0 \le i \le p$, for the corresponding element coming from the $i$'th copy of $V$.
The generators of $E_1^{pq}$ can then be given by the exterior products
$$
\d x^1_{i_1} \wedge \d x^2_{i_2} \wedge \ldots \wedge \d x^q_{i_q},
$$
with $0 \le i_1 \le i_2 \le \ldots \le i_q \le p$.

\begin{nparagraph}
Let us denote by $\pi(n)$ the parity of $n$, i.e., $\pi(n) = 0$ if $n$ is even and $\pi(n) = 1$ if $n$ is odd.

It is immediate to see from the definition of the component $d_1$ of the differential of the complex $TC^\bullet$ that on the first layer of our spectral sequence the differential
can be explicitly written as follows. For $q = 1$ we have
$$
\d x_i \mapsto \pi(i) \d x_{i+1} + (-1)^p \pi(p - i) \d x_i.
$$
For $q = 2$ we have
\begin{align*}
\d x_i \wedge \d y_j \mapsto \pi(i) \d x_{i+1} \wedge \d y_{j+1} + &(-1)^{i+1} \pi(j - i + 1) \d x_i \wedge \d y_{j+1} -\\
&(-1)^{i+1} \delta_{ij} \d y_i \wedge \d x_{i+1} + (-1)^p \pi(p - j) \d x_i \wedge \d y_j.
\end{align*}
Here $i \le j$ and $\delta_{ij}$ is the Kronecker symbol.

Generally, the set of $i_j$ determine an ordered partition of the set of indices $[1, q]$ into $m$ blocks $I_l$, $1 \le l \le m$, such that for any two indices $j$, $k$ within one block
we have $i_j = i_k$, and we will denote this shared value by $v(I_l)$, and moreover if $j \in I_{l_1}$, $k \in I_{l_2}$ with $l_1 < l_2$ then $i_j < i_k$.
Denote by $Sh(i,j)$ the set of $(i,j)$-shuffles, and put by definition $i_0 = -1$ and $i_{q+1} = p + 1$. Then the differential is given by
\begin{multline*}
\d x^1_{i_1} \wedge \d x^2_{i_2} \wedge \ldots \wedge \d x^q_{i_q} \mapsto
\sum_{j=0}^q (-1)^{i_j + 1}\pi(i_{j+1} - i_j + 1) \d x^1_{i_1} \wedge \ldots \wedge \d x^j_{i_j} \wedge \d x^{j+1}_{i_{j+1}+1} \wedge \ldots \wedge \d x^q_{i_q + 1} +\\
\sum_{l=1}^m \quad \sum_{j+k = |I_l| \atop j,k > 0} \  \sum_{\sigma \in Sh(j,k) \atop \sigma \neq \id} (-1)^{v(I_l)+1 + |\sigma|}
\d x^{\sigma(1)}_{i_1} \wedge \ldots \wedge \d x^{\sigma(j)}_{i_j} \wedge \d x^{\sigma(j+1)}_{i_{j+1}+1} \wedge \ldots \wedge \d x^{\sigma(q)}_{i_q + 1}.
\end{multline*}
\end{nparagraph}

For convenience, one could visualize this differential up to the symmetrization summand as edges of a hypercube in a $(q+1)$-dimensional lattice, with elements corresponding to tuples $(p, i_1, \ldots, i_q)$
and the basis vectors corresponding to vectors of the form $(1, 0, \ldots, 0, 1, \ldots, 1)$.

It is clear that generally such hypercubes contribute acyclic pieces of the complex (even when one includes the symmetrization summands), so the only contribution to cohomology
comes from the edge cases when the hypercubes are missing some faces. In fact it only comes from cocycles of the form
$$
\d x^1_0 \wedge \d x^2_1 \wedge \ldots \wedge \d x^{p+1}_p \in E_1^{p,p+1},
$$
which are not fully killed due to the presence of the symmetrization part:
$$
E_2^{p,p+1} = {V^{\tensor (p+1)} \over \Lambda^2 V \tensor V^{\tensor (p-1)} + V \tensor \Lambda^2 V \tensor V^{\tensor (p-2)} + \cdots + V^{\tensor (p-1)} \tensor V} \isom S^{p+1} V.
$$
Putting it all together we find
$$
E_2^{p,q} = \begin{cases}
S^{p+1} V, &\textrm{if $q = p+1$},\\
0, &\textrm{otherwise.}
\end{cases}
$$
The spectral sequence degenerates at the second layer, therefore
$$
TH^\bullet(A) \ \isom\ A \tensor \bigoplus_{p \ge 1} S_A^{p+1}(\Der(A))[-2p-1].
$$
And for the unital $\tensor$-Hochschild complex
$$
UTH^\bullet(A) \ \isom\ A \tensor S_A^\bullet(\Der(A)[-2])[1].
$$
\qed

\begin{nparagraph}[Comparison with the operadic $E_2$-cohomology.]
Let us compare this result to the cohomology of $A$ as an $E_2$-algebra. This can be calculated as the $\Ext$-algebra
$$
\Ext^\bullet_{HH_\bullet(A)}(A, A),
$$
where the Hochschild homology $HH_\bullet(A)$ is considered as an algebra equipped with the shuffle product, and $A$ is a $HH_\bullet(A)$-moduls in a natural way.
In the case of a smooth algebra $A$ the Hochschild-Kostant-Rosenberg theorem gives us isomorphism of algebras
$$
HH_\bullet(A) \ \isom\ \Omega_A^\bullet,
$$
where $\Omega_A^\bullet$ (without de Rham differential) is graded so that $\Omega_A^p$ sits in degree $-p$ and is equipped with the exterior product of differential forms.
The algebra $A$ as $\Omega_A^\bullet$-module is
the quotient $\Omega_A^\bullet / \Omega_A^{\ge 1}$. We find
$$
\Ext^\bullet_{\Omega_A^\bullet}(A, A) \ \isom\ A \tensor S^\bullet(\Der(A)[-2]),
$$
which up to shift by $1$ coincides with our $UTH^\bullet(A)$.

This identification holds in greater generality for any commutative algebra (and even $E_2$-algebra), not just in the smooth case.
\end{nparagraph}

\begin{proposition}
\label{prop_E2_alg}
Let $A$ be a commutative algebra, and $\C = \Perf(A)$ the category of perfect complexes of $A$-modules, then
$$
UTH^\bullet(\C) \ \isom\  UTH^\bullet(A) \ \isom\ HH_{E_2}^\bullet(A)[1] \ \isom\  \Ext^\bullet_{HH_\bullet(A)}(A, A)[1].
$$
\end{proposition}
\proof
Denote by $B$ the reduced Hochschild complex of $A$, equipped with the shuffle product,
$$
B \ =\  \left(\begin{tikzcd}[cramped]
\ldots \ar[r] & A \tensor \wbar{A}^{\tensor 2} \ar[r] & A \tensor \wbar A \ar[r, "0"] & A
\end{tikzcd}
\right).
$$
Under the Dold-Kan equivalence it comes from a simplicial space that we will denote by $sB$ to avoid confusion, and we write $B = DK(sB)$.

The action of $B$ on $A$ is defined by the composition
$$
\begin{tikzcd}[cramped]
B \tensor A \ar[r, epi] & (B / B^{\le -1}) \tensor A \ar[r, "="] & A \tensor A \ar[r, "\mu"] & A.
\end{tikzcd}
$$
The $\Ext$-groups on the right hand side can be expressed using the standard free resolution of $A$ in $B$-modules
$$
P^\bullet \ =\ \left(\begin{tikzcd}[cramped]
\ldots \ar[r] & B^{\tensor_A 3} \tensor_A A \ar[r] & B \tensor_A B \tensor_A A \ar[r] & B \tensor_A A
\end{tikzcd}\right).
$$
Consider the stupid filtration of $P^\bullet$ shifted by one, $F_p(P^\bullet) = \sigma_{\ge(-p-1)} P^\bullet$, and the induced decreasing filtration of the complex $\Hom_B(P^\bullet, A)$.
On the zero layer of the corresponding spectral sequence we have for $p \ge -1$
$$
E_0^{p\bullet} = \Hom_A(B^{\tensor (p+1)}, A).
$$
Denote by $sB^{\boxtimes p}$ the $p$-simplicial space $(\Delta^\op)^p \to \Vect$, defined as
$$
sB^{\boxtimes p}(n_1, \ldots, n_p) = sB(n_1) \tensor \cdots \tensor sB(n_p).
$$
Let $\diag(sB^{\boxtimes p})\from \Delta^\op \to \Vect$ be the totalization of this $p$-simplicial space. Eilenberg-Zilber theorem gives a quasi-isomorphism
between the complex $B^{\tensor p}$ and the $DK(\diag(sB^{\boxtimes p}))$. However, the latter complex is precisely the $p$'th column of the zero layer of the
spectral sequence for the stupid filtration of $UTC^\bullet(A)$ described in the beginning of this section (\ref{equ_sp_seq_alg}). Therefore, the two spectral
sequences for the $\Ext$-groups and $UTH^\bullet(A)$ are isomorphic starting from the first layer, which in turn implies the isomorphism claimed in the proposition.
\qed

\vskip 5em
\subsection{Smooth scheme.}
\label{sec_smooth_scheme}
Let us extend the calculation of the previous section to the case of smooth schemes. First, we establish a simple lemma which will be useful in the future.

\begin{lemma}
\label{lemma_subcomplex}
Assume that $\C$ is generated as a dg-category by an object $X \in \C$, then $TH^\bullet(\C)$ can be calculated using the subcomplex of $TC^\bullet$ generated by $X$:
$$
C^{pq}(\C)_X = \Hom_\k \left( \End_\C(X)^{\tensor (p+1)q}, \End_\C(X^{\tensor (p+1)})\right).
$$
\end{lemma}
\proof
Consider filtration $F$ of the complex $TC^\bullet$ and its restriction to the subcomplex $TC^\bullet_X$. These filtrations induce spectral sequences $E$ and $E_X$ converging to
the cohomology $TH^\bullet(\C)$ and $H^\bullet(TC^\bullet_X)$ respectively. Let us look at the first spectral sequence,
on the layer $E_0$ the column $E_0^{p\bullet}$ is given by the Hochschild complex of $\C^{\tensor(p+1)}$ with
coefficients in the bimodule $\Hom(\boxtimes X_i, \boxtimes Y_i)$. By assumption $\C$ is generated by $X$, in the sense that every object of $\C$, considered as a $\C$-module via Yoneda embedding,
has a resolution, such that every term is of the form $X^{\oplus I}$. It is clear that the tensor product $\C^{\tensor (p+1)}$ is then generated by the exterior product $X^{\boxtimes(p+1)}$.

Therefore, the map of spectral sequences $E_X \to E$, induced by the inclusion $C^{pq}(\C)_X \into C^{pq}(\C)$, is a quasi-isomorphism on the zero layer, hence the
two spectral sequences are isomorphic starting from the first layer. Which immediately implies the statement of the lemma.
\qed

\begin{theorem}
\label{thm_smooth_scheme}
Let $X$ be a smooth scheme, and $\C = D^b(\coh(X))$ with the monoidal structure of derived tensor product of $\O_X$-modules. Then there is a spectral sequence
$$
E_2^{p\bullet} \isom H^\bullet(X, S^{p+1}T_X[-p-1]) \ \Rightarrow\ TH^\bullet(\C),
$$
where $T_X$ is the tangent sheaf of $X$, and $S^n$ denotes the $n$\!'th symmetric power.
In particular, if $\dim X \le 2$, then
$$
TH^\bullet(\C) \ \isom\ \bigoplus_{p \ge 0} H^\bullet(X, S^{p+1} T_X).
$$
\end{theorem}
\proof
Once again, consider filtration $F$ of the complex $TC^\bullet$ and the corresponding spectral sequence. Let $X^{p+1}$ be the Cartesian product of $(p+1)$ copies of $X$.
Since the category $\D^{b}(\coh(X))$ is generated by a single object (see for example \cite{Bondal}), the category $\D^b(\coh(X^{p+1}))$ is generated by external products $X_0 \boxtimes \cdots \boxtimes X_p$ of $(p+1)$ objects of $\C$, then the columns
on the zero layer of the spectral sequence are the Hochschild complexes of $X^{p+1}$. Locally, we are working with a smooth commutative algebra $A$, so as in the
discussion in the beginning of section \ref{subsec_sm_comm} the coefficients are given by $A$, with the structure of $A^{\tensor (p+1)}$-module given by
$$
(a_0 \tensor \cdots \tensor a_p) \cdot b = a_0 \ldots a_p b.
$$
Hence, they glue into the structure sheaf of the diagonal $\Delta\from X \to X^{p+1}$.
Therefore, the first layer of the spectral sequence is
$$
E_1^{p\bullet} \ \isom\ HH^\bullet(X^{p+1}, \O_\Delta).
$$

Since $X$ is smooth, using Hochschild-Kostant-Rosenberg theorem we find
$$
E_1^{p\bullet} \ \isom\ H^\bullet\left(X^{p+1}, \Delta_* \Delta^* (\bigoplus_{q \ge 0} \Lambda^q(T_{X^{p+1}})[-q]) \right).
$$

Locally, the differentials on $E_1$ have already been studied in proposition \ref{prop_sm_comm}. Gluing it together, we obtain
$$
E_2^{p\bullet} \ \isom\ H^\bullet \left( X^{p+1}, \Delta_* (S^{p+1}T_X[-p-1])\right ) \ \isom\ H^\bullet(X, S^{p+1}T_X[-p-1]),
$$
which proves the first statement of the proposition.

Furthermore, if $\dim X \le 2$, then $E_2^{pq}$ are non-zero only for $p+1 \le q \le p+3$, therefore all differentials starting from
$E_2$ vanish and the spectral sequence degenerates.
\qed

Let us give two simple corollaries of this theorem.

\begin{corollary}
\label{cor_sm_scheme}
Let $X$ be a smooth scheme, $\C = D^b(\coh X)$ and $\tensor$ denote the monoidal structure given by the tensor product of coherent sheaves, then the following holds.
\begin{enumerate}[label=\alph*)]
\item
Any infinitesimal automorphism of the underlying category $\C$ preserves the monoidal structure,
and any infinitesimal deformation of $\C$ lifts to a deformation of $(\C, \tensor)$ as monoidal category.
\item
Moreover, there are no infinitesimal automorphisms and deformations of the monoidal structure, preserving the underlying category $\C$.
\end{enumerate}
\end{corollary}
\proof
{\bf a)} The infinitesimal deformations of $\C$ are controlled by the second Hochschild cohomology of $\C$, which is isomorphic to
$H^1(X, T_X)$. This in turn is the component of the second layer of the spectral sequence $E_2^{02}$. Since all differentials from this
component vanish starting from the second layer we find that the canonical map $TH^2(\C) \to HH^2(\C)$ is surjective.

Similarly an infinitesimal automorphisms of $\C$ is given by an element of $E_2^{01} \isom H^0(X, T_X)$ and again all the differentials
starting from $E_2$ vanish. Hence, $TH^1(\C) \to HH^1(\C)$ is surjective.

{\bf b)} It is clear from the theorem that $E_2^{pq} = 0$ for $q \le p$. Therefore,
$$
\wbar{TH}^1(\C) = \wbar{TH}^2(\C) = 0.
$$
\qed

\begin{corollary}
\label{cor_P1}
For $X = \P^1$ and $\C = \D^b(\coh(X))$ we have
$$
TH^{2n+1}(\C) \ \isom\ H^0(\P^1, T_{\P^1}^{\tensor (n+1)}) \ \isom\ H^0(\P^1, \O_{\P^1}(2n + 2)).
$$
\end{corollary}
\proof
Follows from the proposition, since $T_{\P^1} \isom \O(2)$ is an ample line bundle.
\qed

\vskip 5em
\subsection{Associative bialgebra.}
\label{sec_bialgebras}
Let $B$ be an associative and coassociative bialgebra with multiplication $\mu\from B \tensor B \to B$ and coproduct $\Delta\from B \to B \tensor B$.
We denote by $\C$ the derived category of right $B$-modules. Let $M$ and $N$ be two right $B$-modules, the monoidal structure is given by setting
$M \tensor N = M \tensor_\k N$ with the action of $B$ defined as composition
$$
\begin{tikzcd}
(M \tensor N) \tensor B \ar[r, "\Delta"] & (M \tensor N) \tensor (B \tensor B) \ar[r, "\tau_{23}"] & (M \tensor B) \tensor (N \tensor B) \ar[r] & M \tensor N,
\end{tikzcd}
$$
where $\tau_{23}$ denotes the transposition of the second and third factors. Since this tensor product is exact we extend it to $\C$ by taking tensor
products of complexes of right $B$-modules.
Here we implicitly used associativity of the tensor product of vector spaces.
In fact it will be convenient to identify tensor products of vector spaces $(U \tensor V) \tensor W = U \tensor (V \tensor W)$. This way the associators
of the monoidal structure of $\C$ are given by identity maps, due to coassociativity of the coproduct $\Delta$.

Since the associators of $\C$ are identities the higher terms $d_i$ of the differential in the complex $TC^\bullet(\C)$ for $i \ge 2$ all vanish. Therefore,
$TC^\bullet(\C)$ is the totalization of the bicomplex with vertical differential $d_0$ and horizontal differential $d_1$. As before, since $B$ generates
category $\C$ we can restrict our attention to the subcomplex of $TC^\bullet(\C)$ generated by $B$ and we have
$$
C^{pq}(\C)_B = \Hom_\k \left(B^{\tensor (p+1)q}, \End_\C \left(B^{\tensor_\C (p+1)} \right) \right).
$$

\begin{nparagraph}
The deformation of bialgebra $B$ is controlled by the Gerstenhaber-Schack complex $GS^\bullet(B)$, defined as follows. It is the totalization of a bicomplex
$$
GS^{pq} = \Hom_\k(B^{\tensor q}, B^{\tensor p}),
$$
where vertical differentials are the Hochschild differentials induced by multiplication $\mu$ and horizontal differentials are similarly induced by the
comultiplication $\Delta$.

In order to compare the two complexes $TC^\bullet(\C)_B$ and $GS^\bullet(B)$ we introduce a supplementary complex $\wtilde{TC}^\bullet$ to be the totalization
of a tricomplex
$$
\wtilde{TC}^{pqr} = \Hom_\k \left(B^{\tensor (p+1)q}, \Hom_\k \left(B^{\tensor_\C (p+1)} \tensor B^{\tensor r}, B^{\tensor_\C (p+1)}\right) \right).
$$
Here the inside $\Hom$ is the complex expressing $\R\End_{B^\op} \left(B^{\tensor_\C (p+1)} \right)$. Therefore, by construction the natural map $TC^\bullet \to \wtilde{TC}^\bullet$
is a quasi-isomorphism.

Let us rewrite terms $\wtilde{TC}^{pqr}$ in the following way:
$$
\wtilde{TC}^{pqr} = \Hom_\k \left( B^{\tensor r}, \Hom_\k \left(B^{\tensor (p+1)q} \tensor B^{\tensor_\C (p+1)}, B^{\tensor_\C (p+1)}\right) \right).
$$
Now the inside $\Hom$-complex calculates $\R\End_{B^{\tensor (p+1)}} \left( B^{\tensor_\C (p+1)} \right)$, where $B^{\tensor_\C (p+1)}$ is considered as a left $B^{\tensor (p+1)}$-module.
Since it is a free module, the cohomology of this complex is concentrated in degree zero and is isomorphic to $B^{\tensor_\C (p+1)}$.

This way we obtain a map $GS^{p+1,\bullet}(B) \to \wtilde{TC}^{p\bullet\bullet}$, for $p \ge 0$ which is a quasi-isomorphism. Moreover, for $p = -1$ it is immediate to see that
the column $GS^{0\bullet}(B)$ of the Gerstenhaber-Schack complex is quasi-isomorphic to $\R\End_\C(\k)$ and as was discussed in section \ref{sec_unital}, the latter is quasi-isomorphic
to the column $p = -1$ of the unital $\tensor$-Hochschild complex $UTC^\bullet(\C)$. Therefore, the spectral sequences for $GS^\bullet(B)$ and $UTC^\bullet(\C)$ are isomorphic
starting from the first layer (after shift of degree $p$ by 1). And we conclude
$$
UTH^n(\C) \ \isom\ H^{n+1}(GS^\bullet(B)).
$$

\end{nparagraph}

\begin{nparagraph}
In the case when $B$ is a quasi-coassociative bialgebra, the associators of the monoidal structure on $\C$ are no longer identities and the higher terms of the
$\tensor$-Hochschild differential no longer vanish. Our construction can therefore be considered to be a generalization of the Gerstenhaber-Schack complex to the
case of quasi-bialgebras. Using the quasi-isomorphisms of the previous paragraph, it is straightforward to translate these higher terms $d_i$, for $i \ge 2$ from
the language of monoidal categories to the language of quasi-bialgebras.
\end{nparagraph}

%\subsection{Quotient stack $X / G$.}
%Let $X$ be a smooth scheme and $G$ an algebraic group acting on $X$ and consider the quotient stack $X / G$. The tangent complex of $X/G$ has only two non-zero components
%in degrees $-1$ and $0$
%$$
%T_{X/G} = \left(
%\begin{tikzcd}[cramped]
%\g \ar[r, "a"] & T_X
%\end{tikzcd}
%\right).
%$$
%The differential $\g \to T_X$ is given by the infinitesimal action of $\g$ on $X$.

\vskip 5em
\subsection{Kronecker quiver.}
\label{sec_quiver}
Let $Q = (V, E)$ be a finite quiver, where $V$ is the set of vertices and $E$ is the set of arrows. We denote by $s\from E \to V$ the map that picks the beginning of an arrow, and by
$t\from E \to V$ the map that picks the end of an arrow.
Denote by $A = \k Q$ its path algebra, and consider the category of right $A$-modules.
Algebra $A$ contains the commutative subalgebra $\k V$, and is quipped with a bialgebra structure over the enveloping algebra $\k V^e$, with the comultiplication $\Delta \from A \to A \tensor_{\k V^e} A$
given by $\Delta(v) = v \tensor v$ for each vertex $v \in V$ and $\Delta(x) = x \tensor x$ for each arrow $x \in E$. The monoidal structure on the category of right $A$-modules is given by the tensor product
of modules $M \tensor_\C N = M \tensor_{\k V} N$ with the action of $A$ defined by the composition
$$
\begin{tikzcd}
(M \tensor_{\k V} N) \tensor_{\k V} A \ar[r,"\Delta"] &  (M \tensor_{\k V} N) \tensor_{\k V} (A \tensor_{\k V^e} A) \ar[r] & M \tensor_{\k V} N,
\end{tikzcd}
$$
where the second map is given by the action of the first copy of $A$ on $M$ and the second copy of $A$ on $N$.

In other words, a right $A$-module $M$ is a collection of vector spaces $M_v$ indexed by vertices $v \in V$, and linear maps $f^M_x\from M_v \to M_w$ for each arrow $x \from w \to v$. Then the tensor product
of $M$ and $N$ is given by taking tensor product of vector spaces over each vertex $(M \tensor_\C N)_v = M_v \tensor N_v$, and the linear maps
$$
f^{M \tensor_\C N}_x = f^M_x \tensor f^N_x.
$$

Let $\C$ be the bounded derived category of right $A$-modules with the monoidal structure given by the tensor product of complexes of $A$-modules.

Consider a quiver $Q_n$ with two vertices $e_1$ and $e_2$ and $n$ arrows going from $e_2$ to $e_1$.
$$
\begin{tikzpicture}
\def\u{6em}
\node (e1) at (0, 0) {$\bullet$};
\node at (e1) [left] {$e_1$};
\node (e2) at (1*\u, 0) {$\bullet$};
\node at (e2) [right] {$e_2$};
\draw [->] (e2) [out=120,in=60] to node [above] {\scriptsize $x_1$} (e1);
\draw [->] (e2) [out=160,in=20] to node [above] {\scriptsize $x_2$} (e1);
\draw [->] (e2) [out=-120,in=-60] to node [below] {\scriptsize $x_n$} (e1);
\node at (0.5*\u, -0.05*\u) {$\vdots$};
\end{tikzpicture}
$$

\begin{proposition}
\label{prop_quiver}
Let $A = \k Q_n$ for $n \ge 0$, and $\C = A^\op\Mod$, then
$$
TH^\bullet(\C) = 0.
$$
\end{proposition}
\proof
Clearly, $A$ considered as a right module over itself generates category $\C$, so applying lemma \ref{lemma_subcomplex} with $X = A$ we obtain a subcomplex
$$
C^{pq}(\C)_A = \Hom_\k \left( A^{\tensor (p+1)q}, \End_\C\left(A^{\tensor_\C (p+1)}\right)\right),
$$
where the first $A$ is identified with the endomorphism algebra $A \isom \End_\C(A)$ and the tensor product is taken in the category of vector spaces, while the second $A$ is
an object in $\C$ and the tensor product $\tensor_\C$ is taken in the category $\C$.

Modules $e_1 A$ and $e_2 A$ are projective right $A$-modules, generating the category $\C$. We have $A \isom e_1 A \oplus e_2 A$, and
the tensor product of the generators is given by
\begin{align}
\label{equ_tensors}
\begin{split}
e_1 A \tensor e_1 A &\isom e_1 A \oplus e_2 A^{n^2 - n}, \\
e_1 A \tensor e_2 A &\isom e_2 A^n,\\
e_2 A \tensor e_2 A &\isom e_2 A.
\end{split}
\end{align}

Therefore, we have
$$
A^{\tensor 2} = (e_1 A \oplus e_2 A)^{\tensor 2} \isom e_1 A \oplus e_2 A^{n^2 + n + 1},
$$
and generally,
$$
A^{\tensor p} \isom e_1 A \oplus e_2 A^{(n+1)^p - n}.
$$

So the dg-algebra of endomorphisms of $A^{\tensor p}$ is concentrated in degree $0$ and can be written in a block-matrix form
$$
\End_\C(A^{\tensor p}) = \left(
\begin{tikzcd}[sep=small, baseline=-0.25em]
\k & \Hom(\k^{(n+1)^p - n}, \k^n) \\
0 & \gl((n+1)^p - n)
\end{tikzcd}
\right).
$$

First we compute the Hochschild cohomology of $A^{\tensor p}$ with coefficients in $\End_\C(A^{\tensor p})$.

\begin{lemma}
We have
\begin{align*}
HH^0(A^{\tensor p}, \End_\C(A^{\tensor p})) &\isom \k,\\
HH^1(A^{\tensor p}, \End_\C(A^{\tensor p})) &\isom \k^{n^{p+1} - 1},\\
HH^i(A^{\tensor p}, \End_\C(A^{\tensor p})) &= 0, \quad \textit{for $i > 0$}.
\end{align*}
\end{lemma}
\proof
Since $A$ is a quiver algebra, we have the standard resolution of $A$ by projective $A$-bimodules of length one
$$
\begin{tikzcd}
\displaystyle
\bigoplus_{x_i} A e_1 \tensor \k x_i \tensor e_2 A \arrow[r, mono] & (A e_1 \tensor e_1 A) \oplus (A e_2 \tensor e_2 A) \ar[r, epi] & A.
\end{tikzcd}
$$

Therefore, we obtain a projective resolution $K^\bullet$ of the tensor product $A^{\tensor p}$ in $A^{\tensor p}$-bimodules of length $p$. Let $(z_1, \ldots, z_p)$
be an ordered set of indices, where each $z_j$ is either $e_1$, $e_2$ or one of $x_i$. We put $|z_j| = 0$ if $z_j = e_k$, $|z_j| = 1$ if $z_j = x_i$ and write
$|z_1, \ldots, z_p| = \sum_j |z_j|$ for the number of $x_i$'s in the collection.
For such a collection we consider a left $A$-module
$$
L(z_1, \ldots, z_p) = A t(z_1) \boxtimes \cdots \boxtimes A t(z_p),
$$
and a right $A$-module
$$
R(z_1, \ldots, z_p) = s(z_1) A \boxtimes \cdots \boxtimes s(z_p) A.
$$
Here, by convention we put $s(e_i) = t(e_i) = e_i$.

So the component in degree $k$ of the resolution of $A^{\tensor p}$ is given by
\begin{equation}
\label{equ_K_complex}
K^{-k} \ \isom\ \bigoplus_{(z_1, \ldots z_p) \atop |z_1, \ldots, z_p| = k} L(z_1, \ldots, z_p) \tensor \k(z_1, \ldots z_p) \tensor R(z_1, \ldots, z_p),
\end{equation}
where $\k(z_1, \ldots z_p)$ is a one-dimensional vector space spanned by the symbol $(z_1, \ldots z_p)$.

Intuitively, one can visualize this as a $p$-dimensional cube, with $k$-dimensional faces indexed by collections $(z_1, \ldots z_p)$ with
$|z_1, \ldots z_p| = k$. So each one-dimensional edge is decorated by $n$ arrows, $2$-dimensional face by $n^2$ diagonal arrows, etc.

Now, let us look at the coefficients $E = \End_\C(A^{\tensor p})$. It is a $(\k V)^{\tensor p}$-bimodule, so it can also be visualized as arrows
in the $p$-dimensional cube. Moreover, since we are interested in the space of $A^{\tensor p}$-bimodule maps from $K^\bullet$ to $E$ we can restrict
our attention only to arrows going in the same direction as those in $K^\bullet$, as described above. Let us describe this explicitly.

Consider $A^{\tensor p}$ as a $(A^{\tensor p}, A)$-bimodule. The left action makes it into a $V^p$-graded right $A$-module.
For a vertex of the cube labeled by a collection $e_\bullet = (e_{i_1}, \ldots e_{i_p})$, containing at least one $e_2$, the subspace of
$A^{\tensor p}$ over this vertex is isomorphic to the direct sum of $e_2 A$ with the basis given by symbols $[z_1, \ldots, z_p]$, such that $z_k = e_2$ if $e_{i_k} = e_2$ and
$z_k$ is one of $x_i$ if $e_{i_k} = e_1$. We will denote this subspace $W(e_\bullet)$. For the remaining vertex $(e_1, \ldots, e_1)$ the corresponding $A$-module is isomorphic to
\begin{equation}
\label{equ_W'}
e_1 A \oplus \bigoplus_{[z_1, \ldots, z_p]} e_2 A_{[z_1, \ldots, z_p]},
\end{equation}
where the sum is taken over all collections $[z_1, \ldots, z_p]$, such that each $z_k$ is one of $x_i$'s, and $z_k$'s are not all equal to each other.
We will denote this sum over $[z_1, \ldots, z_p]$ by $W'$.
In the case when $z_1 = \ldots = z_p = x_i$ we will write $e_2 A[z_1, \ldots, z_p]$ for the corresponding submodule of $e_1 A$ isomorphic to $e_2 A$.
And denote the sum
$$
W'' = \bigoplus_{z_1 = \ldots = z_p = x_i} e_2 A_{[z_1, \ldots, z_p]}.
$$

The space of endomorphisms $E$ is graded in the following way, over a vertex $(e_{i_1}, \ldots e_{i_p})$ containing at least one $e_2$, we have
$\End(W(e_\bullet))$, and over $(e_1, \ldots, e_1)$ the block-matrix algebra
$$
\begin{pmatrix}\k & \Hom(W', W'') \\ 0 & \End(W') \end{pmatrix}.
$$
Over an arrow from $e_\bullet$ to $e'_\bullet$, when neither vertex is $(e_1, \ldots, e_1)$ we have $\Hom(W(e_\bullet), W(e'_\bullet))$, and finally,
from $(e_1, \ldots, e_1)$ to $e_\bullet$
$$
\Hom(W', W(e_\bullet)),
$$
and from $e_\bullet$ to $(e_1, \ldots, e_1)$
$$
\Hom(W(e_\bullet), W' \oplus W'').
$$

\begin{nparagraph}
Let us illustrate this in the case of $p = 2$ and $n = 2$. For the complex $K^\bullet$ we have
\begin{equation}
\label{equ_K_example}
\begin{tikzcd}[sep=4em]
(e_1, e_2) \ar[d, shift left=0.25em] \ar[d, shift right=0.25em] & (e_2, e_2) \ar[l, shift left=0.25em, "{(x_2, e_2)}"] \ar[l, shift right=0.25em, "{(x_1, e_2)}"'] \ar[d, shift left=0.25em] \ar[d, shift right=0.25em] \ar[ld, Rightarrow, "4"] \\
(e_1, e_1) & (e_2, e_1). \ar[l, shift left=0.25em] \ar[l, shift right=0.25em]
\end{tikzcd}
\end{equation}
The vertices are in cohomological degree $0$, horizontal and vertical arrows are in degree $-1$ and the four diagonal arrows are in degree $-2$.

The subspaces $W$ of $A^{\tensor 2}$ introduced above are
\begin{align*}
W(e_2, e_2) &= \k_{[e_2, e_2]},\\
W(e_2, e_1) &= \k_{[e_2, x_1]} \oplus \k_{[e_2, x_2]},\\
W(e_1, e_2) &= \k_{[x_1, e_2]} \oplus \k_{[x_2, e_2]},\\
W' &= \k_{[x_1, x_2]} \oplus \k_{[x_2, x_1]},\\
W'' &= \k_{[x_1, x_1]} \oplus \k_{[x_2, x_2]}.\\
\end{align*}

To simplify the picture we will only draw the part of $E$ corresponding to arrows going in the same direction as in the picture (\ref{equ_K_example}).
$$
\begin{tikzcd}[column sep=8em, row sep=4em]
{\End(W(e_1, e_2))} \ar[d, "{\Hom(W(e_1, e_2), W' \oplus W'')}"'] & {\End(W(e_2, e_2))} \ar[l, "{\Hom(W(e_2, e_2), W(e_1, e_2))}"']  \ar[d, "{\Hom(W(e_2, e_2), W(e_2, e_1))}"] \ar[ld, sloped, "{\Hom(W(e_2, e_2), W' \oplus W'')}"] \\
{\begin{pmatrix}\k & \Hom(W', W'') \\ 0 & \End(W') \end{pmatrix}} & {\End(W(e_2, e_1))}.  \ar[l, "{\Hom(W(e_2, e_1), W' \oplus W'')}"]
\end{tikzcd}
$$
\end{nparagraph}

\begin{nparagraph}
Take a Hochschild cochain $\phi \from K^{-k} \to E$, then the value of $d\phi$ on an arrow $z_\bullet$ with $|z_\bullet| = k + 1$ is given by
$$
d\phi(z_\bullet) = \sum_{j, |z_j| = 1} (-1)^{\epsilon_j} \left( z_j \mathop{\circ}\limits_j \phi(\ldots, z_{j-1}, e_2, z_{j+1}, \ldots) - \phi(\ldots, z_{j-1}, e_1, z_{j+1}, \ldots) \mathop{\circ}\limits_j z_j \right),
$$
where the sum is taken over all $z_j$ in the collection $z_\bullet$ that are equal to some $x_i$, and $\epsilon_j = |z_1, \ldots, z_{j-1}|$. The action of $x_i$ on the righthand side
comes from the $A^{\tensor p}$-bimodule structure on $E$, induced by the left action on $A^{\tensor p}$. Explicitly,
$$
x_i \mathop{\circ}\limits_j [z_1, \ldots, z_{j-1}, e_2, z_{j+1}, \ldots z_p] = [z_1, \ldots, z_{j-1}, x_i, z_{j+1}, \ldots z_p].
$$

Consider a filtration $F$ of the Hochschild complex $C^\bullet = \Hom(K^\bullet, E)$ by the distance of an arrow from the vertex $(e_2, \ldots, e_2)$. In other words, we put
$$
\delta(z_1, \ldots, z_p) = |z_1, \ldots, z_p| + \#\left\{ z_j \mid z_j = e_1 \right\}.
$$
We have an increasing filtration $F_m K^\bullet$ formed by restricting the sum over $(z_1, \ldots, z_p)$ in (\ref{equ_K_complex}) to arrows with $\delta(z_1, \ldots, z_p) < m$.
It induces a dual decreasing filtration
$$
F_m C^\bullet = \Ker(C^\bullet \to \Hom(F_m K^\bullet, E)).
$$

The part $\gr^0_F C^\bullet$ is isomorphic to $\k$ positioned in cohomological degree $0$. From the previous description of the cochains and the differential, it is clear that the associated graded pieces
$\gr^i_F(C^\bullet)$ are acyclic for $1 \le i < p$. And the contribution of $\gr^p_F(C^\bullet)$ is concentrated in degree one and is isomorphic to the quotient
\begin{equation}
\label{equ_cohom_class}
{ \Hom\left(\oplus_{x_i} \k(e_1, \ldots, e_1, x_i), \Hom(W(e_1, \ldots, e_1, e_2), W' \oplus W'')\right) \over \k \oplus \Hom(W', W'') \oplus \End(W')}.
\end{equation}
Here $\k \oplus \Hom(W', W'') \oplus \End(W')$ is a subspace of $\End(W' \oplus W'')$, and the $\Hom$-space at the top is isomorphic to
$$
\Hom \left( \oplus_{x_i} \k(e_1, \ldots, e_1, x_i) \tensor W(e_1, \ldots, e_1, e_2), W' \oplus W'' \right) \isom \End(W' \oplus W'').
$$
Therefore,
$$
\dim HH^1(A^{\tensor p}, E) = n^{2p} - ( (n^p - n)^2 + n(n^p - n) + 1 ) = n^{p+1} - 1.
$$
\end{nparagraph}
\qed

\begin{nparagraph}
The lemma says that the first layer of the spectral sequence computing the $TH^\bullet(\C)$ has only two non-zero rows.
$$
E_1^{p0} \isom \k, \quad E_1^{p1} \isom \k^{n^{p+2} - 1}, \quad p \ge 0.
$$
The differential $d^p\from E_1^{p0} \to E_1^{p+1, 0}$ on the zero row is an alternating sum of $p+3$ terms, each equal to identity map on $\k$, therefore it is an isomorphism if $p$ is even, and zero if $p$ is odd. Hence the cohomology
in the zero row all vanish.

Now let us look at the first row. First let us rewrite the quotient in (\ref{equ_cohom_class}) as
$$
\Hom(W'', W' \oplus W'') / \k \cdot \Id_{W''}.
$$
The basis of the $\Hom$ space  is formed by elements $[x_{i_1}, \ldots, x_{i_p}]_j$, which correspond to the map sending
$$
(e_1, \ldots, e_1, x_j) \mapsto \left( [x_j, \ldots, x_j, e_2] \mapsto [x_{i_1}, \ldots, x_{i_p}] \right).
$$
In this notation the identity on $W''$ is written as
$$
\sum_j [x_j, \ldots, x_j]_j.
$$
Pick a cocycle representative of the cohomological class $[x_{i_1}, \ldots, x_{i_p}]_j$ of the following form
$$
(e_1, \ldots, x_j, \ldots, e_1) \mapsto \left( [x_j, \ldots, e_2, \ldots x_j] \mapsto [x_{i_1}, \ldots, x_{i_p}] \right),
$$
where $x_j$ in the first symbol and $e_2$ in the second symbol are placed in the same position $1 \le k \le p$. In other words, we send all arrows ending at the vertex $(e_1, \ldots, e_1)$ and marked by $x_j$ to the same element
$[x_{i_1}, \ldots, x_{i_p}]$, and send all other arrows to zero.
\end{nparagraph}

\begin{nparagraph}
Let us explicitly describe the differential on $E_1$. For a cocycle
$$
\phi \in E_1^{p-1,1} = HH^1(A^{\tensor p}, \End_\C(A^{\tensor p})),
$$
we have
\begin{multline*}
d\phi(a_1, \ldots, a_{p+1}) = a_1 \tensor \phi(a_2, \ldots, a_{p+1}) - \phi((a_1 \tensor a_2), \ldots, a_{p+1}) + \cdots +\\ (-1)^{p} \phi(a_1, \ldots, (a_p \tensor a_{p+1})) + (-1)^{p+1} \phi(a_1, \ldots, a_p) \tensor a_{p+1}.
\end{multline*}
Let us clarify notation here. Since the right $A$-module $A$ generates the category $\C$ the Hochschild cochain $\wtilde\phi$, representing class $\phi$ extends to a Hochschild cochain $\wtilde\phi^+$ of the category $\C^{\tensor p}$,
uniquely up to coboundary. We consider the tensor product $(a_1 \tensor a_2) \tensor \cdots \tensor a_{p+1}$ as an element of $\End_\C(A \tensor A) \tensor \End_\C(A)^{\tensor (p-1)}$,
and apply to it the extended cochain $\wtilde\phi^+$. Moreover, since the tensor product $A \tensor A$ splits into a direct sum of direct summands of $A$, construction of such an extension $\wtilde\phi^+$ is
especially simple. Indeed, consider a projective module $P = e_1 A \tensor_\k M \oplus e_2 A \tensor_\k N$, for some vector spaces $M$ and $N$. Then an endomorphism of $P$ can be written as
$$
f = f_1 e_1 + f_2 e_2 + \sum_{i} g_i x_i,
$$
where $f_1 \in \End(M)$, $f_2 \in \End(N)$ and $g_i \in \Hom(N, M)$. If $\psi \in C^1(A, A)$ is a Hochschild cochain, then the extension
$$
\wtilde\psi^+(f) = f_1 \psi(e_1) + f_2 \psi(e_2) + \sum_{i} g_i \psi(x_i).
$$

We will apply this to an arrow $z_\bullet = (e_{i_1}, \ldots, e_{i_k-1}, x_j, e_{i_{k+1}}, \ldots, e_{i_{p+1}})$ and the cocycle representing $\phi = [x_{i_1}, \ldots, x_{i_p}]_j$, as described above.
Then we use formulas (\ref{equ_tensors}) to decompose the tensor product $A \tensor A$ into a direct sums of indecomposable projective modules $e_k A$.

Specifically, since $e_2 A \tensor e_2 A \isom \bigoplus e_2 A$, we have $\phi(a_1, \ldots, (e_2 \tensor e_2), \ldots, a_{p+1})$ is a sum of $\phi(a_1, \ldots, e_2, \ldots, a_{p+1})$. However, our cocycle $\phi$ vanishes whenever there is at least one $e_2$ among it's arguments.
Therefore, $d\phi$ vanishes whenever there are at least two $e_2$'s among its arguments.

Similarly, for $e_2 \tensor x_j$ and $x_j \tensor e_2$ we use the previous decomposition and $e_2 A \tensor e_1 A \isom \bigoplus e_2 A$. So the map once again decomposes into sum of $e_2$'s and the cocycle $\phi$ vanishes. The same applies to $e_2 \tensor e_1$ and $e_1 \tensor e_2$,
therefore, if $z_\bullet$ contains exactly one $e_2$ then $d\phi(z_\bullet)$ vanishes when $e_2$ is neither in the first place nor the last, and
\begin{align*}
d\phi(z_\bullet) &= e_2 \tensor \phi(z_2, \ldots, z_{p+1}), \quad \textrm{if $z_1 = e_2$},\\
d\phi(z_\bullet) &= (-1)^{p+1} \phi(z_1, \ldots, z_p) \tensor e_2, \quad \textrm{if $z_{p+1} = e_2$}.
\end{align*}

For $e_1 \tensor e_1$ we use decomposition $e_1 A \tensor e_1 A \isom e_1 A \oplus \bigoplus e_2 A$. So the map decomposes into a sum of a single $e_1$ and several $e_2$. Again, since cocycle $\phi$ vanishes if there is at least one $e_2$ the only surviving term can
be obtained by replacing $e_1 \tensor e_1$ with $e_1$. In terms of $d\phi$ this implies that if the arrow $z_\bullet = (e_1, \ldots, x_j, \ldots, e_1)$ with $x_j$ in position $l$, then terms $\phi(z_1, \ldots, (z_k \tensor z_{k+1}), \ldots z_{p+1})$ with $k \neq (l-1)$ and $k \neq l$
in the formula for the differential evaluated at $[x_j, \ldots, e_2, \ldots, x_j]$, where $e_2$ is in position $l$, equal to $[x_{i_1}, \ldots, x_{i_{k-1}}, x_{i_k}, x_{i_k}, x_{i_{k+1}}, \ldots, x_{i_p}]$.

Finally for $e_1 \tensor x_j$ (and similarly $x_j \tensor e_1$) we have a map
$$
e_1 A \tensor e_2 A \isom \bigoplus_i e_2 A_{x_i} \to e_1 A \oplus W' \isom e_1 A \tensor e_1 A,
$$
where $W'$ was defined in (\ref{equ_W'}). Arguing as before, cocycle $\phi$ vanishes on the component landing in $W'$, so we only need to consider part landing in $e_1 A$. Recall that submodules of $e_1 A \subset A \tensor A$ isomorphic to $e_2 A$ are indexed by
symbols $[x_i, x_i]$. Therefore, the only contribution to $\phi$ comes from the map $x_j\from e_2 A_{x_j} \to e_1 A$. So in terms of the differential $d\phi$, if $z_\bullet$ is as before terms $\phi(z_1, \ldots, (z_k \tensor z_{k+1}), \ldots z_{p+1})$ with $k = (l-1)$ or $k = l$
evaluated at $[x_j, \ldots, e_2, \ldots, x_j]$ are again equal to $[x_{i_1}, \ldots, x_{i_{k-1}}, x_{i_k}, x_{i_k}, x_{i_{k+1}}, \ldots, x_{i_p}]$.
\end{nparagraph}

\begin{nparagraph}
\label{par_quiver_E1}
To describe the element $d\phi$ in the cohomology group $HH^1(A^{\tensor p}, \End_\C(A^{\tensor p}))$, using the calculation of the cohomology from the previous lemma, we observe that
it is enough to evaluate $d\phi$ on two arrows: $(e_1, \ldots, e_1, x_j)$ and $(x_j, e_1, \ldots, e_1, e_2)$ and add their contributions. The contribution to $d\phi$ coming from the first arrow is
$$
[x_j, x_{i_1}, \ldots, x_{i_p}]_j - [x_{i_1}, x_{i_1}, x_{i_2}, \ldots, x_{i_p}]_j + \cdots + (-1)^p [x_{i_1}, \ldots, x_{i_p}, x_{i_p}]_j.
$$
And the contribution from the second arrow is
$$
(-1)^{p+1} [x_{i_1}, \ldots, x_{i_p}, x_j]_j.
$$

The acyclicity of the complex $E_1^{p\bullet}$ then follows from the lemma.
\end{nparagraph}

\begin{lemma}
Consider the semi-free associative dg-algebra on $n$ generators $R_n = \k\!\<x_1, \ldots, x_n\>$, with $\deg x_i = 1$, and differential given by $dx_i = x_i^2$. Then
$$
H^\bullet(R) = \k.
$$
\end{lemma}
\proof
First of all, notice that if $n = 1$, then the dg-algebra $R_1$ is $(\k\!\<x\>, dx = x^2)$ and the statement is obvious. In general, the dg-algebra $R_n$ is isomorphic to
the coproduct in the category of associative dg-algebras of $n$ copies of $R_1$ and the statement follows from the following simple observation.

Let $P$ and $Q$ be two associative dg-algebras concentrated in non-negative degrees, such that $H^\bullet(P) \isom \k$ and $H^\bullet(Q) \isom \k$, furthermore, assume that
$Q$ is augmented, and denote by $\wbar Q$ the augmentation ideal. Let $P * Q$ be their coproduct, then we have decomposition
$$
P * Q = P \oplus \left (P \tensor \wbar Q \tensor P \right) \oplus \bigoplus_{k \ge 2} P \tensor \left( \wbar Q \tensor P\right)^{\tensor k}.
$$
By assumption, the first summand is quasi-isomorphic to $\k$. Since $\wbar Q$ is acyclic, using Kunneth formula we conclude that the rest of the summands are also acyclic, hence
$P * Q$ is quasi-isomorphic to $\k$.
\qed

Now, consider the complex described in paragraph \ref{par_quiver_E1}. Let us embed it into $R_n$ by mapping $E_1^{p1} \to R_n^{p+3}$ as follows:
$$
[x_{i_1}, \ldots, x_{i_p}]_j \mapsto x_j x_{i_1} x_{i_2} \ldots x_{i_p} x_j.
$$
Clearly, this embedding is compatible with differentials on $E_1$ and in $R_n$. Moreover, the image of the embedding is a direct summand of complex $R_n$, hence  $E_1^{p1}$ is acyclic,
as $R_n$ is acyclic in degrees greater or equal than three. Therefore, the spectral sequence degenerates starting from $E_2$ and all the cohomology groups vanish.

This completes the proof of proposition \ref{prop_quiver}.
\qed

\begin{nparagraph}
It may be of a particular interest to look at the case of the quiver $Q_2$. In this case the bounded derived category of right $A$-modules is equivalent to the bounded derived category of coherent sheaves on the projective line $\P^1$,
by sending module $e_2 A$ to the structure sheaf $\O_{\P^1}$ and $e_1 A$ to the twisted line bundle $\O_{\P^1}(1)$. To avoid confusion we will write $\C_{Q_2}$ for the category $\C$ with the monoidal structure given by the tensor product of representations of the quiver.
From the proof of proposition \ref{prop_quiver} we see that the cohomology group $\wbar{TH}^1(\C_{Q_2})$ is trivial and $\wbar{TH}^2(\C_{Q_2})$
is three-dimensional. Moreover, since both $TH^1(\C_{Q_2})$ and $TH^2(\C_{Q_2})$ vanish from the long exact sequence of cohomology groups
$$
\begin{tikzcd}
\ldots \ar[r] & TH^1(\C_{Q_2}) \ar[r] & HH^1(\C_{Q_2}) \ar[r] & \wbar{TH}^2(\C_{Q_2}) \ar[r] & TH^2(\C_{Q_2}) \ar[r] & \ldots
\end{tikzcd}
$$
we see that the natural map $HH^1(\C_{Q_2}) \to \wbar{TH}^2(\C_{Q_2})$ is an isomorphism. This map describes the action of the Lie algebra of infinitesimal automorphisms of $\C$ on the space of monoidal structures,
therefore all infinitesimal deformations of the monoidal structure on $\C$ are obtained by applying infinitesimal automorphisms of the
underlying category.

Let us compare this to the case of the monoidal structure on $\C$ given by the tensor product of coherent sheaves discussed in section \ref{sec_smooth_scheme}. We will denote this monoidal category by $\C_{\P^1}$. Combining results of corollaries \ref{cor_sm_scheme} and \ref{cor_P1} we have
$\wbar{TH}^1(\C_{\P^1}) = \wbar{TH}^2(\C_{\P^1}) = 0$, $TH^2(\C_{\P^1}) = 0$ and $TH^2(\C_{\P^1})$ is three-dimensional. So the natural map $HH^1(\C_{\P^1}) \to \wbar{TH}^2(\C_{\P^1})$ is zero, in other words this monoidal structure is rigid and the Lie algebra of
infinitesimal automorphisms of $\P^1$ acts trivially on it.

In particular we find that these two monoidal structures belong to two different orbits of $\Aut(\C)$ acting on the space of monoidal structures on $D^b(\coh \P^1)$.

\end{nparagraph}

\vfill\eject


\begin{thebibliography}{WWW}
\addcontentsline{toc}{section}{References}
\bibitem[BB]{Bondal} A. Bondal, M. van den Bergh {\em Generators and representability of functors in commutative and noncommutative geometry.} Moscow Mathematical Journal {\bf 3} (2003) Number 1, pp. 1--36.
\bibitem[CL]{ChuLaz} J. Chuang, A. Lazarev {\em L-infinity maps and Twistings.} Homology, Homotopy and Applications, {\bf 13(2)} (2011) pp. 175--195.
\bibitem[CY]{Yetter} L. Crane, D. N. Yetter. {\em Deformations of (bi)tensor categories.} Cahiers de Topologie et Géométrie Différentielle Catégoriques (1998).
\bibitem[Dav]{Davydov} A. Davydov. {\em Twisting of monoidal structures.} \href{http://arxiv.org/abs/q-alg/9703001}{http://arxiv.org/abs/q-alg/9703001}
\bibitem[Ger]{Gerstenhaber} M. Gerstenhaber. {\em The Cohomology Structure of an Associative Ring.} Annals of Mathematics, {\bf 78}, No. 2 (1963), pp. 267--288.
\bibitem[GHS]{Gainutdinov} A. M. Gainutdinov, J. Haferkamp, C. Schweigert. {\em  Davydov-Yetter cohomology, comonads and Ocneanu rigidity.} Advances in Mathematics {\bf 414} (2023).
\bibitem[GS]{Gerst-Schack} M. Gerstenhaber, S. D. Schack. {\em Algebras, bialgebras, quantum groups, and algebraic deformations.} Contemp. Math. {\bf 134}, Amer. Math. Soc. (1992), pp. 51--92.
\bibitem[GY]{Ginot} G. Ginot, Sinan Yalin. {\em Deformation theory of bialgebras, higher Hochschild cohomology and formality.} \href{https://arxiv.org/abs/1606.01504}{https://arxiv.org/abs/1606.01504}.
\bibitem[KL]{Keller} B. Keller, W. Lowen. {\em On Hochschild Cohomology and Morita Deformations.}  International Mathematics Research Notices (2009) No. 17, pp. 3221--3235.
\bibitem[KS]{KoSo} M. Kontsevich, Y. Soibelman. {\em Deformation Theory. I} \href{http://www.math.ksu.edu/~soibel/Book-vol1.ps}{http://www.math.ksu.edu/~soibel/Book-vol1.ps}
\bibitem[PS]{Panero} P. Panero, B. Shoikhet. {\em The category $\Theta_2$, derived modifications, and deformation theory of monoidal categories.} \href{https://arxiv.org/abs/2210.01664v2}{https://arxiv.org/abs/2210.01664v2}
\bibitem[Shr]{Shrestha} T. Shrestha. {\em Algebraic deformations of monoidal category.} Kansas State University, PhD thesis (2010).
\bibitem[Sta]{Stasheff} J. Stasheff. {\em Drinfeld's quasi-Hopf algebras and beyond.} Contemp. Math. {\bf 134}, Amer. Math. Soc. (1992), pp. 297--307.
\bibitem[Yet]{Yetter2} D. N. Yetter. {\em Braided deformations of monoidal categories and {V}assiliev invariants.} Higher category theory ({E}vanston, {IL}, 1997),  117--134, Contemp. Math., 230, Amer. Math. Soc., Providence, RI, 1998.
\end{thebibliography}
\end{document}